\documentclass[letterpaper,11pt]{amsart}
\usepackage[active]{srcltx}
\usepackage{longtable}
\usepackage{amsmath,amssymb,amscd, booktabs}
\usepackage[all]{xy}
\usepackage[headings]{fullpage}
\newcommand{\A}{{\mathbb A}}
\newcommand{\Q}{{\mathbb Q}}
\newcommand{\Z}{{\mathbb Z}}
\newcommand{\R}{{\mathbb R}}
\newcommand{\C}{{\mathbb C}}

\newcommand{\p}{\mathfrak p}
\newcommand{\OF}{{\mathfrak o}}
\newcommand{\GL}{{\rm GL}}
\newcommand{\PGL}{{\rm PGL}}

\newcommand{\SO}{{\rm SO}}

\newcommand{\GSp}{{\rm GSp}}

\newcommand{\St}{{\rm St}}

\newcommand{\val}{{\rm val}}

\newcommand{\Norm}{{\rm N}}
\newcommand{\Mat}{{\rm M}} 
\newcommand{\Gal}{{\rm Gal}}
\newcommand{\GSO}{{\rm GSO}}
\newcommand{\GO}{{\rm GO}}

\newcommand{\Ind}{{\rm Ind}}
\newcommand{\Trace}{{\rm Tr}}

\newcommand{\SSp}{{\rm Sp}}

\newcommand{\mat}[4]{{\setlength{\arraycolsep}{0.5mm}\left[
\begin{array}{cc}#1&#2\\#3&#4\end{array}\right]}}
\newcommand{\forget}[1]{}
\newcommand{\nl}{

\vspace{2ex}}
\newcommand{\nll}{

\vspace{1ex}}
\def\qdots{\mathinner{\mkern1mu\raise0pt\vbox{\kern7pt\hbox{.}}\mkern2mu
\raise3.4pt\hbox{.}\mkern2mu\raise7pt\hbox{.}\mkern1mu}}

\allowdisplaybreaks

\newtheorem{lemma}{Lemma.}[section]

\newtheorem{proposition}[lemma]{Proposition}

\newtheorem*{maintheorem}{Main Theorem}

\begin{document}

\title[Siegel modular forms from Hilbert modular forms]{Siegel modular forms of degree two \\ attached to Hilbert modular forms}
\author[Johnson-Leung and Roberts]{Jennifer Johnson-Leung\\
Brooks Roberts}
\address{Department of Mathematics, University of Idaho\\Moscow, Idaho, 83843 USA}
\begin{abstract}
Let $E/\Q$ be a real quadratic field and $\pi_0$ a cuspidal, irreducible, automorphic representation of $\GL(2,\A_E)$ with trivial central character and infinity type $(2,2n+2)$ for some non-negative integer $n$.  We show that there exists a non-zero Siegel paramodular newform $F: \mathfrak H_2 \to \mathbb C$ with weight, level, Hecke eigenvalues, epsilon factor and $L$-function determined explicitly by $\pi_0$.  We tabulate these invariants in terms of those of $\pi_0$ for every prime $p$ of $\Q$.
\end{abstract}
\maketitle

\section{Introduction}
In his 1980 paper \cite{Y}, Yoshida studied certain explicit theta liftings of Hilbert modular forms of weight $(2, 2n+2)$ for real quadratic extensions of $\Q$ to Siegel modular forms of degree 2 and weight $n+2$ for the Siegel congruence subgroup $\Gamma_0(N)$ and an appropriate Dirichlet character $\chi$.  Yoshida calculated the action of the Hecke operators $T(1,1,p,p)$ and $T(1,p,p,p^2)$, defined below, on these lifts for $p\nmid N$, though Yoshida did not determine when these lifts are non-zero.
\nl
In this paper, we study an analogous problem.  Given a Hilbert modular form of weight $(2, 2n+2)$ we prove the existence of a non-zero Siegel modular form of degree 2 and weight $n+2$ for the paramodular congruence subgroup.  Our main theorem completely characterizes the resulting Siegel modular form, including the Hecke eigenvalues at {\em every} rational prime $p$.
\begin{maintheorem}
Let $E$ be a real quadratic extension of $\mathbb Q$ with real archimedean places $\infty_1$ and $\infty_2$. Let $\pi_0$ be a cuspidal irreducible automorphic representation of $\GL(2,\A_E)$ with trivial central character. Let $\mathfrak{N}_0$ be the level of $\pi_0$. Assume that $\pi_0$ is not Galois-invariant and $\pi_{0,\infty_1} = D_2$ and $\pi_{0,\infty_2} = D_{2n+2}$ with $n \geq 0$ a non-negative integer and $D_k$ the holomorphic discrete series representation of  $\PGL(2,\R)$ of lowest weight $k$. Let $N=d_{E}^2\Norm_{\mathbb Q}^E(\mathfrak{N}_0)$, where $d_E$ is the discriminant of $E/\mathbb Q$. Then there exists a non-zero  Siegel paramodular newform $F: \mathfrak H_2 \to \mathbb C$ of weight $k=n+2$ and paramodular level $N$ such that:
\begin{enumerate} 
\item For every prime $p$, 
\begin{equation}
\label{heckeFeq}
T(1,1,p,p)F =p^{k-3} \lambda_p  F \quad \text{and} \quad T(1,p,p,p^2)F= p^{2(k-3)}  \mu_p F
\end{equation}
where the Hecke eigenvalues $\lambda_p$ and $\mu_p$ are determined by the Hecke eigenvalues of $\pi_0$ as follows.  If $p$ splits, let $w_1$ and $w_2$ be the places above $p$.  If $p$ is non-split, let $w$ be the place above $p$.
\begin{enumerate}
\item If $\val_p(N)=0$, 
$$
\lambda_p=
\begin{cases}
p(\lambda_{w_1} +\lambda_{w_2})& \text{if $p$ is split,}\\
0&\text{if $p$ is not split,}
\end{cases}
\quad
\mu_p=
\begin{cases}
p^2+p\lambda_{w_1}\lambda_{w_2}-1 &\text{if $p$ is split,}\\
-(p^2+p\lambda_w  +1) & \text{if $p$ is not split}. 
\end{cases}
$$
\item If $\val_p(N)=1$, then $p$ splits and $\val_{w_1}(\mathfrak{N}_0)=1$, 
$\val_{w_2}(\mathfrak{N}_0)=0$, and
$$
\lambda_p=p\lambda_{w_1} + (p+1)\lambda_{w_2}, \qquad \mu_p=p\lambda_{w_1}\lambda_{w_2}.
$$
\item If $\val_p(N) \geq 2$, then:
\begin{enumerate}
\item[] $p$ inert:
$$
\lambda_p=0, \qquad \mu_p=-p^2-p\lambda_w;
$$ 
\item[] $p$ ramified:
$$
\lambda_p=p \lambda_w, \qquad 
\mu_p=
\begin{cases}
0 & \text{if $\val_w(\mathfrak{N}_0)=0$},\\
-p^2 & \text{if $\val_w(\mathfrak{N}_0) \geq 1$};
\end{cases}
$$
\item[] $p$ split and $\val_{w_1}(\mathfrak{N}_0) \leq \val_{w_2}(\mathfrak{N}_0)$:
$$
\lambda_p=p(\lambda_{w_1}+\lambda_{w_2}),\qquad 
\mu_p=
\begin{cases}
0&\text{if $\val_{w_1}(\mathfrak{N}_0) =0$,}\\
-p^2&\text{if $\val_{w_1}(\mathfrak{N}_0)\geq 1$}.
\end{cases}
$$
\end{enumerate}
\end{enumerate}
For particular $\pi_0$, 
 $\lambda_p$ and $\mu_p$ are given  in Proposition \ref{Hecketable}. 
\item For every prime $p |N$, let $U_p$ be the Atkin-Lehner operator, defined below.  Then, 
\begin{equation}
\label{ALFeq}
 F|_k U_p = \varepsilon_p F
\end{equation}
with 
$$
\varepsilon_p =
\begin{cases}
 \varepsilon(1/2,\pi_{0,w_1},\psi_p,dx_{\psi_p}) \varepsilon(1/2,\pi_{0,w_2},\psi_p,dx_{\psi_p})& \mbox{if $p$ is split,}\\
\varepsilon(1/2,  \pi_{0,w}, \psi_w,dx_{\psi_w})\omega_{E_w/\mathbb{Q}_p}(-1)&\mbox{if $p$ is not split,}
\end{cases}
$$
where $\psi_w$ is an additive character of $E_w$ with conductor $\OF_w$.  For particular $\pi_0$, $\varepsilon_p$ is given in Proposition \ref{Hecketable}. 
\item For every prime $p$, we have an equality of Euler factors 
\begin{equation}\label{peulerfactor}
L_p(s+k-\frac{3}2, F) = L_p(s,\pi_0), 
\end{equation}
where $k=n+2$ and $L_p(s,F)$ is defined below for every finite place $p$ of $\Q$.  Moreover, we have the functional equation
\begin{equation}\label{functionaleq}
\Lambda(2k-2-s,F)=(-1)^{k}(\prod_{p\mid N}\varepsilon_p)N^{s-k+1}\Lambda(s,F)
\end{equation}
where the completed $L$-function is defined as the product
$$\Lambda(s,F)=(2\pi)^{-2s}\Gamma(s)\Gamma(s-k+2)\prod_{p<\infty}L_p(s,F).$$
\end{enumerate}
\end{maintheorem}
Our main theorem has potential applications to arithmetic geometry.  For example, Consani and Scholten \cite{CS} studied the four dimensional Galois representation $\rho$ of $G_\Q$ on the \'etale cohomology $H^3(\tilde{X}_{\bar{\Q}},\Q_\ell(\sqrt{5}))$ of the desingularization $\tilde{X}$ of a quintic three-fold $X$. Consani and Scholten showed that $\rho$ is induced from a representation $\sigma$ of $\Gal(\bar \Q/ \Q(\sqrt{5}))$.  By work of Yi \cite{Yi}, $\sigma$ corresponds to an automorphic representation $\pi_0$, as in our main theorem, with $n=1$ and $\mathfrak{N}_0=\p_2\p_3\p_5^2=(30)$.  As a consequence of our main theorem, there exists a non-zero Siegel modular form $F$ of degree 2 and weight 3 for the paramodular group of level $N=2^23^25^4$ such that $\Lambda(s,F)=\Lambda(s,\rho)$.
\nl
The paper is organized as follows.  In section 2, we introduce notation and definitions which will be used throughout the paper.  Section 3 supplies the proof of the main theorem.  This proof depends on certain local results which are proved in sections 4 and 5.  In section 4, we explicitly tabulate local the $L$-packets, Hecke eigenvalues, and epsilon factors used in the proof of the main theorem.  The technical heart of the paper is in section 5, where we calculate the gamma factors of the Novodvorsky zeta integrals of a generic supercuspidal representation of $\GSp(4,F)$ for a non-archimedean local field $F$. 

\section{Notation and Definitions} 
Let 
$$
J= \begin{bmatrix} & \mathbf{1}_2 \\
-\mathbf{1}_2 &
\end{bmatrix}.
$$
We define the algebraic $\mathbb Q$-group $\GSp(4)$ as the set of all $g \in \GL(4)$ such that ${}^t g J g = x J$ for some $x \in \GL(1)$; we call $x$ the multiplier of $g$ and denote it by $\lambda (g)$. The kernel of $\lambda: \GSp(4) \to \GL(1)$ is the symplectic group $\SSp(4)$. For $N$ a positive integer we define 
$$
\Gamma^{\mathrm{para}} (N) =
\begin{bmatrix}
\mathbb Z & N \mathbb Z & \mathbb Z &\mathbb Z \\
\mathbb Z & \mathbb Z & \mathbb Z & N^{-1} \mathbb Z \\
\mathbb Z & N \mathbb Z & \mathbb Z & \mathbb Z \\
N \mathbb Z & N \mathbb Z & N \mathbb Z & \mathbb Z
\end{bmatrix} \cap \SSp(4,\mathbb Q).
$$
Let $p$ be a prime of $\mathbb Q$. For $n \geq 0$ a non-negative integer define the local paramodular group $\mathrm{K}(p^n)$ as the group of all $k \in \GSp(4,\mathbb Q_p)$ such that 
$$
k \in
\begin{bmatrix}
\mathbb Z_p & p^n \mathbb Z_p  & \mathbb Z_p &\mathbb Z_p \\
\mathbb Z_p & \mathbb Z_p & \mathbb Z_p & p^{-n} \mathbb Z_p \\
\mathbb Z_p & p^n \mathbb Z_p & \mathbb Z_p & \mathbb Z_p \\
p^n \mathbb Z_p & p^n \mathbb Z_p & p^n \mathbb Z_p & \mathbb Z_p
\end{bmatrix} 
$$
and $\det (k) \in \mathbb Z_p^\times$. Note that $\mathrm{K}(p^0)= \GSp(4,\mathbb Z_p)$. We have
$$
\Gamma^{\mathrm{para}} (N) = \GSp(4,\mathbb Q) \cap
\GSp(4,\mathbb R)^+ \prod_{p < \infty} \mathrm{K}(p^{\val_p(N)}).
$$
Here, $p^{\val_p(N)}$ is the exact power of $p$ dividing $N$, and  $\GSp(4,\mathbb R)^+$ is the subgroup of $g \in \GSp(4,\mathbb R)$ such that $\lambda (g) > 0$. Let $\mathfrak H_2$ be the Siegel upper half-space of degree two. Then $\GSp(4,\mathbb R)^+$ acts on $\mathfrak H_2$ by 
$$
g  \langle Z \rangle
=
(AZ+B)(CZ+D)^{-1}, \qquad g= \begin{bmatrix} A & B \\ C & D \end{bmatrix}, \quad Z \in \mathfrak H_2
$$
and we define $j(g,Z)=\det (CZ+D)$. If $k$ is a positive integer, $g \in \GSp(4,\mathbb R)^+$,  and $F:\mathfrak H_2 \to \mathbb C$ is a  function we define
$$
(F|_kg)(Z) = \lambda(g)^k j(g,Z)^{-k} F(g\langle Z \rangle), \quad
Z \in \mathfrak H_2.
$$
A {\it Siegel modular form} of degree $2$ and weight $k$ with respect to $\Gamma^{\mathrm{para}}(N)$ is a holomorphic function $F:\mathfrak H_2 \to \mathbb C$ such that $F|_k\gamma =F$ for $\gamma \in \Gamma^{\mathrm{para}}(N)$.  Let $M_k(\Gamma^{\mathrm{para}}(N))$ and $S_k(\Gamma^{\mathrm{para}} (N))$ be the spaces of all Siegel modular forms or cuspforms of degree $2$ and weight $k$ with respect to $\Gamma^{\mathrm{para}}(N)$, respectively.  For each prime $p$, we define two Hecke operators $T(1,1,p,p)$ and $T(1,p,p,p^2)$ on $M_k(\Gamma^{\mathrm{para}}(N))$  as follows. Let 
$$
\Gamma^{\mathrm{para}}(N) 
\begin{bmatrix}
1&&&\\
&1&&\\
&&p&\\
&&&p
\end{bmatrix}
\Gamma^{\mathrm{para}} (N) = \bigsqcup \Gamma^{\mathrm{para}} (N) 
h_i
$$
and
$$
\Gamma^{\mathrm{para}}(N) 
\begin{bmatrix}
p&&&\\
&1&&\\
&&p&\\
&&&p^2
\end{bmatrix}
\Gamma^{\mathrm{para}} (N) = \bigsqcup \Gamma^{\mathrm{para}} (N) 
h'_j
$$
be disjoint decompositions. 
Note that if $p \nmid N $ then  
$$
\Gamma^{\mathrm{para}}(N)\mathrm{diag}(p,1,p,p^2) \Gamma^{\mathrm{para}}(N)=\Gamma^{\mathrm{para}}(N)\mathrm{diag}(1,p,p^2,p) \Gamma^{\mathrm{para}}(N).
$$
For $F \in M_k(N)$ set
$$
T(1,1,p,p) F= p^{k-3}\sum_i F|_k h_i,\quad T(1,p,p,p^2) F= p^{2(k-3)}\sum_j F|_k h_j'.
$$
If $N=1$, this definition is the same as in, for example, (1.3.3) of \cite{A1}.
For each prime $p$ dividing $N$, choose a matrix $\gamma_p \in \SSp(4,\mathbb Z)$ such that
$$
\gamma_p \equiv 
\begin{bmatrix}
&&&1\\
&&1&\\
&-1&&\\
-1&&&
\end{bmatrix}
\quad \mathrm{mod}\ p^{\val_p(N)} \quad \mathrm{and} \quad \gamma_p 
\equiv
\begin{bmatrix}
1&&&\\
&1&&\\
&&1&\\
&&&1
\end{bmatrix} \quad \mathrm{mod} \ Np^{-\val_p(N)},
$$
and define
$$
U_p=\gamma_p 
\begin{bmatrix}
p^{\val_p(N)}&&&\\
&p^{\val_p(N)}&&\\
&&1&\\
&&&1
\end{bmatrix}.
$$
It can be verified that $U_p$ normalizes $\Gamma^{\mathrm{para}}(N)$ and that $U_p^2$ is contained in $p^{\val_p(N)} \Gamma^{\mathrm{para}}(N)$, so that $F \mapsto F|_k U_p$ defines an involution of $S_k(\Gamma^{\mathrm{para}}(N))$ and $M_k(\Gamma^{\mathrm{para}}(N))$.  Let $S_k^{\mathrm{new}}(\Gamma^{\mathrm{para}}(N))$ be defined as in \cite{RS2}. Let $F \in S_k^{\mathrm{new}}(\Gamma^{\mathrm{para}}(N))$ and assume that
$$
T(1,1,p,p)F =\lambda_{F,p} F, \quad T(1,p,p,p^2)F=\mu_{F,p} F, \quad F|_k U_p = \varepsilon_{F,p} F. 
$$
Then we define $L_p(s,F)$ as follows:
\begin{enumerate}
\item If $\val_p(N)=0$, then
$$
L_p(s, F)^{-1}=1-\lambda_{F,p}p^{-s}+(p\mu_{F,p}+p^{2k-3}+p^{2k-5})p^{-2s}-p^{2k-3}\lambda_{F,p}p^{-3s}+p^{4k-6}p^{-4s}.
$$  
\item If $\val_p(N)=1$, then
$$
L_p(s,F)^{-1}=1-(\lambda_{F,p}+p^{k-3}\varepsilon_{F,p})p^{-s}+(p\mu_{F,p}+p^{2k-3})p^{-2s}+\varepsilon_{F,p}p^{3k-5}p^{-3s}.
$$
\item If $\val_p(N)\geq2$, then
$$
L_p(s,F)^{-1}=1-\lambda_{F,p}p^{-s}+(p\mu_{F,p}+p^{2k-3})p^{-2s}.
$$
\end{enumerate}
Note that the case $\val_p(N)=0$  agrees with the classical Euler factor at $p$ as given in \cite{A1}, Theorem 3.1.1.  The work \cite{A1} uses $T(p^2)$ instead of $T(1,p,p,p^2)$. However, one has the relation $pT(1,p,p,p^2)+p^2(p+1)T(p,p,p,p)=T(p)^2-T(p^2)-p^2T(p,p,p,p)$ by 3.3.38 of \cite{A2}. Compare also Theorem 2 of \cite{Sh}. The definitions in the cases $\val_p(N) \geq 1$ are motivated by the results of \cite{RS1}. 
\nl
{\bf Additional notation:} For $k$ a positive integer, we let $D_k$ denote the holomorphic discrete series representation of $\PGL(2,\R)$ of lowest weight $k$. Suppose that $L$ is a nonarchimedean local field of characteristic zero with ring of integers $\OF$ and prime ideal $\p=\varpi \OF \subset \OF$. We define the character $\nu:L^\times \to \C^\times$ by $\nu(x)=|x|$ where $|\cdot|$ is the absolute value such that $|\varpi|=|\OF/\p|^{-1}$. If $\chi:L^\times \to \C^\times$ is a character, then  $a(\chi)$ is the smallest non-negative integer $n$ such that $\chi(1+\p^n)=1$, where we take $1+\p^0 = \OF^\times$. Let $(\tau,V)$ be a generic, irreducible, admissible representation of $\GL(2,L)$ with trivial central character. For $n$ a non-negative integer, let $\Gamma_0(\p^n)$ be the subgroup of $\left[ \begin{smallmatrix} a&b \\ c&d \end{smallmatrix} \right]$ in $\GL(2,\OF)$ such that $c \equiv 0\ (\p^n)$. 
The group $\Gamma_0(\p^n)$ is normalized by the Atkin-Lehner element $\left[\begin{smallmatrix} &1 \\ -\varpi^n& \end{smallmatrix} \right]$. 
We define $a(\tau)$ to be the smallest non-negative integer $n$ such that $V^{\Gamma_0(\p^n)} \neq 0$; we call $\p^{a(\tau)}$ the level of $\tau$. The space $V^{\Gamma_0(\p^{a(\tau)})}$ is spanned by a non-zero vector $v$. We have $\tau (\left[\begin{smallmatrix} &1 \\ -\varpi^n& \end{smallmatrix} \right])v = \varepsilon_\tau v$ for $\varepsilon_\tau = \varepsilon(1/2,\tau,\psi) \in \{\pm 1\}$, where $\psi: L \to \C$ is a  character with conductor $\OF$. We call $\varepsilon_\tau$ the Atkin-Lehner eigenvalue of $\tau$. We also have $T(\p) v = \lambda_\tau v$ for some $\lambda_\tau \in \C$. Here, $T(\p)v = \sum_i \tau(h_i)v$, where $\Gamma_0(\p^{a(\tau)}) \left[\begin{smallmatrix} \varpi&\\ & 1 \end{smallmatrix} \right] \Gamma_0(\p^{a(\tau)}) = \sqcup h_i \Gamma_0(\p^{a(\tau)})$ is a disjoint decomposition. We call $\lambda_\tau$ the Hecke eigenvalue of $\tau$. The group $\GSp(2n)$ is defined with respect to  $\left[\begin{smallmatrix} &1_n\\ -1_n \end{smallmatrix}\right]$, and $1_{\GSp(2n)}$ and $\St_{\GSp(2n)}$ are the trivial and Steinberg representations of $\GSp(2n,L)$, respectively.

\section{Proof of the Main Theorem}

{\bf Proof of the Main Theorem:} We begin by constructing a certain cuspidal, irreducible, admissible representation $\pi$ of $\GSp(4,\A)$ with trivial
central character; the desired Siegel modular form will correspond to a particular vector in $\pi$. 
For every place $v$ of $\mathbb Q$, define $\pi_{0,v}=\bigotimes_{w|v}\pi_{0,w}$.
Let $\varphi(\pi_{0,v}):W_{F_v}'\rightarrow \GSp(4,\mathbb C)$ be the $L$-parameter associated to $\pi_{0.v}$ as in (\ref{spliteq}) and (\ref{nonspliteq}); if $v$ is non-split in $E$, we take $\eta=1$. 
By \cite{B} the representation $\pi_{0,w}$ is tempered for all finite places $w$ of $E$. Let $\Pi(\varphi(\pi_{0,v}))$ be the $L$-packet of tempered, irreducible, admissible representations of $\GSp(4,\Q_v)$ with trivial central character associated to $\varphi(\pi_{0,v})$ as in \cite{R}.  For finite $v=p$, the packet $\Pi(\varphi(\pi_{0,p}))$ coincides with the packet associated to $\varphi(\pi_{0,p})$ in \cite{GT}, and
contains a unique generic representation $\pi_p$ of $\GSp(4,\mathbb Q_p)$.  It is known that $\Pi(\varphi(\pi_{0,\infty}))$ contains the lowest weight representation $\pi_k$ with $k=n+2$; for the precise definition of $\pi_k$, see \cite{AS}, p. 184. We set $\pi_\infty=\pi_k$.  Define $\pi$ to be the restricted tensor product
$$
\pi=\bigotimes_v\pi_v.
$$
By Theorem 8.6 of \cite{R}, $\pi$ is a cuspidal, irreducible, admissible, automorphic representation of $\GSp(4,\A_\Q)$ with trivial central character. 
\nl
To define the appropriate vector in $\pi$,  for each finite prime $p$ of $\mathbb Q$, let $\Phi_p$ be the local paramodular newform in $\pi_p$, and let $\Phi_\infty \in \pi_\infty$ be the non-zero smooth vector as in Lemma 3.4.2 of \cite{AS}. 
Note that $\Phi_p$ exists and is unique up to scalars by \cite{RS1}; we may assume that for almost all $p$, $\Phi_p$ is the unramified vector used to define the restricted tensor product. Also, $\pi_\infty( u)\Phi_\infty =j(u,I)^{-k} \Phi_\infty$ for $u \in U(2)$, where
$$
I = \mat{i}{}{}{i} \in \mathfrak H_2
$$
and $U(2)$ is the subgroup of $u \in \SSp(4,\R)$ such that $u\langle I \rangle = I$. 
We set
$$
\Phi = \bigotimes_v \Phi_v.
$$
For each finite prime $p$ of $\mathbb Q$, because $\Phi_p$ is a local paramodular newform in $\pi_p$,  we have $T_{1,0}(p)\Phi_p = \lambda_p \Phi_p$ and $T_{0,1}(p)\Phi_p= \mu_p \Phi_p$ for some 
complex numbers $\lambda_p$ and $\mu_p$ and $\pi_p(u_p) \Phi_p =\varepsilon_p \Phi_p$ for some $\varepsilon_p \in \{ \pm 1\}$; here, $T_{0,1}(p)$ and $T_{1,0}(p)$ are
the Hecke operators from Chapter 6 of \cite{RS1} and $u_p$ is the Atkin-Lehner element defined in (2.2) of \cite{RS1}. In fact, 
$\lambda_p$, $\mu_p$ and $\varepsilon_p$ are as in i) and ii)  by Proposition \ref{Hecketable}.
\nl
Next, define $F: \mathfrak H_2 \to \mathbb C$ by $F(Z)=\lambda (h)^{-k}j(h,I)^{k} \Phi (h_\infty)$ where 
$h \in \GSp(4,\mathbb R)^+$ is such that $h\langle I \rangle = Z$. 
Then $F$ is 
holomorphic by Lemma 3.2.1 of \cite{AS}, and an argument shows that $F \in S_{k}^{\mathrm{new}}( \Gamma^{\mathrm{para}}(N))$. 
 A computation  shows that
$$
T(1,1,p,p)F = \lambda_p p^{k-3} F, \qquad T(1,p,p,p^2)F= \mu_p p^{2(k-3)} F.
$$
A similar computation shows that $F|_k U_p =\varepsilon_p F$ because 
$\pi_p(u_p) \Phi_p = \varepsilon_p \Phi_p$. This proves i) and ii). To prove iii), we note that the equality (\ref{peulerfactor}) follows by comparing the Euler factors
at each finite prime $p$ of $\mathbb Q$ and using (\ref{heckeFeq}), (\ref{ALFeq}) and (\ref{Lequalitieseq}). To deduce (\ref{functionaleq}), we recall the functional equation for the completed $L$-function of $\pi_0$ (e.g., Theorem 6.2 of \cite{J}),
$$
L(1-s,\pi_0)=\varepsilon(1-s,\pi_0)L(s,\pi_0).
$$
For every rational prime $p$, we have a canonical additive character $\psi_p(x)=e^{-2\pi i\lambda(x)}$ where $\lambda$ is the composition
$$
\lambda:\Q_p\rightarrow\Q_p/\Z_p\hookrightarrow\Q/\Z\hookrightarrow\R/\Z,
$$
and $\psi_\infty(x)=e^{2\pi ix}$.  Then the function $\psi(x)=\prod_v\psi_v(x_v)$ defines a character of $\A/\Q$.  For the quadratic extension $E$, we  have the character $\tilde{\psi}=\psi\circ\Trace_{E/Q}$ of $\A_E/E$.  For each finite prime, $w$ of $E$, we define $n(\tilde{\psi_w})$ to be the least non-negative integer $n$ such that $\tilde{\psi_w}(\p_w^{-n})=1$.  The epsilon factor is computed via
\begin{align*}
\varepsilon(s,\pi_0)&=(-1)^n\prod_{w<\infty}\varepsilon(s,\pi_{0,w},\tilde{\psi_w })\quad\text{(\cite{Ge}, Theorem 6.16)}\\
&=(-1)^n\prod_{p|N}\prod_{w|p}\varepsilon(\varphi_{0,w},\tilde{\psi_w},dx_{\tilde{\psi_w}})q_w^{-s(2n(\tilde{\psi_w})+a(\varphi_{0,w}))}\quad\text{(\cite{Rohr} 11 Prop.)}\\
&=(-1)^n\prod_{\substack{p|N\\\text{split}}}\varepsilon(\varphi_p,\psi_p,dx_{\psi_p})p^{-sa(\varphi_p)}\prod_{\substack{p|N\\ \text{ nonsplit}}}\varepsilon(\varphi_{0,w},\tilde{\psi_w},dx_{\tilde{\psi_w}})q_w^{-s(2n(\tilde{\psi_w})+a(\varphi_{0,w}))}\\
&=(-1)^n\prod_{\substack{p|N\\\text{split}}}\varepsilon(\varphi_p,\psi_p,dx_{\psi_p})p^{-sa(\varphi_p)}\\&\qquad\times\prod_{\substack{p|N\\\text{ nonsplit}}}\varepsilon(\varphi_{0,w},\tilde{\psi_w},dx_{\tilde{\psi_w}})q_w^{-2sn(\tilde{\psi_w})}p^{-s(a(\varphi_p)-2d(E_w/\Q_p))}\quad\text{(by (\ref{Npieq}))}\\
&=(-1)^n\prod_{p|N}\varepsilon(\varphi_p,\psi_p,dx_{\psi_p})p^{-sa(\varphi_p)}\quad\text{(\cite{Rohr} 11 Prop., $q_w^{n(\tilde{\psi_w})}=p^{d(E_w/\Q_p)}$)}\\
&=(-1)^n\prod_{p|N}\varepsilon_pp^{-a(\varphi_p)(s-1/2)}\quad\text{(by (\ref{epsiloneq}))}\\
&=(-1)^nN^{1/2-s}\prod_{p|N}\varepsilon_p.
\end{align*}
Finally, the archimedean Euler factors of $\pi_0$ are given by
$$
L(s,\pi_{0,\infty_1})L(s,\pi_{0,\infty_2})=(2\pi)^{-2s-n-1}\Gamma(s+1/2)\Gamma(s+(2n+1)/2).
$$
Substituting into the functional equation for $\pi_0$ yields (\ref{functionaleq}).
\qed

\section{Local Results}\label{local}

\subsection*{Some definitions} Throughout this section let $F$ be a  nonarchimedean local field of characteristic zero with ring of
integers $\OF$, let $\p$ be the maximal ideal of  $\OF$ with $\p= \varpi \OF$, and let $q$ be the number
of elements of $\OF/\p$. Also, let $E$ be a quadratic extension of $F$ with $\Gal(E/F)=\{1,\sigma\}$ and associated quadratic character $\omega_{E/F}$. The residue class degree of $E/F$ is denoted by $f=f(E/F)$ and the valuation of the discriminant by $d=d(E/F)$.  Let
$$
J'=
\begin{bmatrix}
&&&1\\
&&1&\\
&-1&&\\
-1&&&
\end{bmatrix},
\qquad
K=\begin{bmatrix}
1&&&\\
&1&&\\
&&&1\\
&&1&
\end{bmatrix}.
$$
We note that in the work \cite{RS1} the group $\GSp(4)$ is defined with respect $J'$, while throughout this work
we define $\GSp(4)$ with respect to $J$, as in section 2. However, it is easy to see that conjugation by $K$ defines
an isomorphism between the two realizations, in either direction.

\subsection*{Two families of $L$-parameters for $\GSp(4)$}
We now consider two families of $L$-parameters for $\GSp(4)$ over $F$. The
 first family, which we will call the split case,  is parameterized by pairs $(\pi_1,\pi_2)$ where $\pi_1$
and $\pi_2$ are irreducible, admissible representations
of $\GL(2,F)$ having the same central character, while the second family, which
we refer to as the non-split case,  is parameterized by triples 
$(E, \pi_0,\eta)$ where $E$ is a quadratic extension of $F$, 
$\pi_0$ is an irreducible, admissible representation of $\GL(2,E)$
with Galois invariant central character $\omega_{\pi_0}$, and $\eta$ is a character
of $F^\times$ such that $\omega_{\pi_0} = \eta \circ \Norm_F^E$. 
\nl
To define the parameter $\varphi(\pi_1,\pi_2):W_F' \to \GSp(4,\mathbb C)$ associated to a pair $(\pi_1,\pi_2)$,
let $\varphi_1: W_F' \to \GL(2,\mathbb C)$ and $\varphi_2: W_F' \to \GL(2,\mathbb C)$
be the $L$-parameters of $\pi_1$ and $\pi_2$, respectively. We define
\begin{equation}
\varphi(\pi_1,\pi_2) (x)
=
\begin{bmatrix}
a_1&&b_1&\\
&a_2&&b_2\\
c_1&&d_1&\\
&c_2&&d_2
\end{bmatrix} \quad  \text{for}\ 
\varphi_1(x) = \mat{a_1}{b_1}{c_1}{d_1}, \
\varphi_2(x) = \mat{a_2}{b_2}{c_2}{d_2} \ \text{and} \ x \in W_F'. \label{spliteq}
\end{equation}
Thus, $\varphi(\pi_1,\pi_2)$ is the symplectic direct sum of $\varphi_1$ and $\varphi_2$. 
\nll
To define the $L$-parameter $\varphi(E,\pi_0,\eta):W_F' \to \GSp(4,\mathbb C)$
associated to a triple $(E,\pi_0,\eta)$,
let $\varphi_0:W_E' \to \GL(2,\mathbb C)$ be the $L$-parameter of $\pi_0$, and let $V_0$ be the space of $\varphi_0$.  Let $g_0$ be a representative for the nontrivial coset of $W'_E\backslash W'_F$.  We consider the representation of $W'_F$ induced from $\varphi_0$ which we can realize as $V_0\oplus V_0$ via the isomorphism
$$
\Ind_{W'_E}^{W'_F} \varphi_0 \stackrel{\sim}{\longrightarrow} V:=V_0 \oplus V_0
$$
that sends $f$ to $f(1) \oplus f(g_0)$. 

Notice that $\eta$ can also be considered as a character of $W'_F$ with the property that $\eta|_{W'_E}=\det(\varphi_0)$.

We define a non-degenerate symplectic bilinear form on $V$ by
$$
\langle v_1 \oplus v_2, v_1' \oplus v_2' \rangle
=
\eta(g_0) \langle v_1, v_1' \rangle + \langle v_2,v_2' \rangle.
$$
A computation shows that $\langle x\cdot v,x\cdot v' \rangle = \eta(x) \langle v,v' \rangle$
for $v,v'\in V$ and $x \in W'_F$.  The $L$-parameter $\varphi(E,\pi_0,\eta)$ is defined to be the representation $V$ with this symplectic structure.
We can choose a symplectic basis for $V$ so that for $y\in W'_E$, 
\begin{gather}
\varphi(E,\pi_0,\eta) (y)
=
\begin{bmatrix}
a&&\eta(g_0)^{-1}b&\\
&a'&&b'\\
\eta(g_0)c&&d&\\
&c'&&d'
\end{bmatrix} \quad  \text{for}\ 
\varphi_0(y) = \mat{a}{b}{c}{d} \ \text{and}\ 
\varphi_0(g_0yg_0^{-1}) = \mat{a'}{b'}{c'}{d'},\nonumber \\
\text{and} \quad
\varphi(E,\pi_0,\eta) (g_0)
=
\begin{bmatrix}
&1&&\\
a_0&&\eta(g_0)^{-1}b_0&\\
&&&\eta(g_0)\\
c_0&&\eta(g_0)^{-1}d_0&
\end{bmatrix} \quad  \text{for}\ 
\varphi_0(g^2_0) = \mat{a_0}{b_0}{c_0}{d_0}. \label{nonspliteq}
\end{gather}

\subsection*{Associated $L$-packets.}
The local Langlands conjecture is proven for $\GSp(4)$ over $F$ in \cite{GT}, and in the following proposition we tabulate the $L$-packets $\Pi(\varphi)$ associated to the $L$-parameters $\varphi$ constructed above for all choices of $\pi_1$, $\pi_2$ and $\pi_0$. 
To list these $L$-packets we proceed as follows. First, the work \cite{RS1} gives an explicit map, determined
by the desiderata of the local Langlands conjecture for $\GSp(4,F)$, from the set of non-supercuspidal,
irreducible, admissible representations of $\GSp(4,F)$ to the set of $L$-parameters; moreover, this map is the same as that
in \cite{GT}. Using this map, it is straightforward to determine all the non-supercuspidal elements in 
the $L$-packets $\Pi(\varphi)$ for $\varphi$ as above. We remark that in some cases, to use \cite{RS1}, it is necessary to consider an $L$-parameter equivalent to $\varphi$. 
In \cite{RS1} the non-supercuspidal representations are divided into eleven groups based on inducing 
data.  Note that in particular each group contains a generic representation, which is designated by the letter ``a'' if the group contains more than one type of irreducible representation.  We also use  the notation of \cite{ST} and \cite{RS1} for the representations in Table 1, below.
\nl
 Next, the $L$-packets $\Pi(\varphi)$ containing supercuspidal
elements for $\varphi$ as above can be described as follows. Assume that $\pi_1$ and $\pi_2$ are discrete series representations such that $\pi_1 \ncong \pi_2$. Then the $L$-packet $\Pi(\varphi(\pi_1,\pi_2))$ consists of two elements. These two representations of $\GSp(4,F)$ are theta lifts 
of the form $\theta_{X}(\sigma)$ and $\theta_{X'}(\sigma')$ where $X$ and $X'$ are the
four dimensional hyperbolic and anisotropic quadratic spaces, respectively, and $\sigma$ and $\sigma'$ are irreducible, admissible representations
of $\GO(X,F)$ and $\GO(X',F)$, respectively, arising from the pair $\pi_1$ and $\pi_2$. It is known that $\theta_X(\sigma)$ and $\theta_{X'}(\sigma')$ have central character $\omega_{\pi_1}=\omega_{\pi_2}$,  $\theta_X(\sigma)$ is generic and tempered, and $\theta_{X'}(\sigma)$ is non-generic. When both $\pi_1$ and $\pi_2$ are not supercuspidal, i.e., $\pi_1= \alpha \St_{\GL(2)}$ and $\pi_2=\beta\St_{\GL(2)}$ for some characters $\alpha$ and $\beta$ of $F^\times$, then $\theta_X(\sigma)=\delta([\alpha^{-1}\beta,\nu\alpha^{-1}\beta],\nu^{-1/2}\alpha)$; this representation belongs to group Va. 
On the other hand, $\theta_{X'}(\sigma')$ is supercuspidal, and to supplement the partition from \cite{RS1}, we will say that it belongs to group $\text{Vb}^*$. If exactly one of $\pi_1$ and $\pi_2$ is supercuspidal, say $\pi_1=\alpha \St_{\GL(2)}$ with $\alpha$ a character of $F^\times$, then  $\theta_X(\sigma) =\delta(\nu^{1/2} \alpha^{-1} \pi_2,\nu^{-1/2}\alpha)$; this representation is of type XIa. Again, $\theta_{X'}(\sigma')$ is supercuspidal, and we say that it is of type $\text{XIb}^*$.  If both $\pi_1$ and $\pi_2$ are supercuspidal, then both $\theta_X(\sigma)$ and $\theta_{X'}(\sigma')$ are supercuspidal, and we say that they are of type $\text{XIIa}^*$ and $\text{XIIb}^*$, respectively. 
Finally, assume that $(E,\pi_0,\eta)$ is as above with $\pi_0$ supercuspidal and not Galois invariant, i.e., $\pi_0^{\sigma} \ncong \pi_0$. Then the $L$-packet $\Pi(\varphi(E,\pi_0,\eta))$ consists of a single representation. This representation is a theta lift of the form $\theta_{X_0}(\sigma_0)$ where $X_0$ is the four dimensional quadratic space having anisotropic component $(E,\Norm_F^E)$ and $\sigma_0$ is the supercuspidal, irreducible, admissible representation of $\GO(X_0,F)$ associated to $\pi_0$ and $\eta$. The representation $\theta_{X_0}(\sigma_0)$ is generic, supercuspidal and has central character $\eta$; we say that this representation is in group $\text{XIII}^*$.  Note that the superscript ${}^*$ indicates that representations of a given type are supercuspidal.
\nl
The table of the following proposition summarize the data obtained by this method.   
\begin{proposition} 
\label{Lpackettable}
Let $\varphi=\varphi(\pi_1, \pi_2)$ (resp, $\varphi(E, \pi_0, \eta)$) be as defined in the previous subsection.  Then the $L$-packet $\Pi(\varphi)$ associated to $\varphi$ by the local Langlands correspondence is given in Table I (resp. II) below.

\end{proposition}
\begin{center}
\renewcommand{\arraystretch}{1.5}
\begin{longtable}{lll}
\caption{$L$-Packets}\\
\toprule
condition & $\Pi(\varphi)$&Group \\
\hline
\endfirsthead

\multicolumn{3}{c}{ \tablename\ \thetable{} -- continued from previous page}\\
\hline
condition & $\Pi(\varphi)$& Group \\
\endhead

\multicolumn{3}{c}{ Continued on next page}
\endfoot

\bottomrule
\endlastfoot

\multicolumn{3}{l}{\textbullet\qquad\qquad
$\pi_1 = \alpha_1 \times \alpha_2$\qquad $\pi_2=\alpha_1' \times \alpha_2'$}\\
\hline
$\alpha_1'\alpha_2^{-1}\neq \nu^{\pm 1}$, $\alpha_2' \alpha_2^{-1} \neq \nu^{\pm 1}$ &
$\alpha_1' \alpha_2^{-1}\times \alpha_2'\alpha_2^{-1} \rtimes \alpha_2$& I \\
\hline
$\alpha_1'\alpha_2^{-1}\neq \nu^{\pm 1}$, $\alpha_2' \alpha_2^{-1} = \nu$ &
$\alpha_1'\alpha_2^{-1} \rtimes \nu^{1/2}\alpha_21_{\GSp(2)}$&IIIb \\
\hline
$\alpha_1'\alpha_2^{-1}\neq \nu^{\pm 1}$, $\alpha_2' \alpha_2^{-1} = \nu^{-1}$ &
$\alpha_2'\alpha_1^{-1} \rtimes \nu^{1/2}\alpha_1 1_{\GSp(2)}$&IIIb \\
\hline
$\alpha_1'\alpha_2^{-1}= \nu$,  $\alpha_2' \alpha_2^{-1} \neq \nu^{\pm 1}$ &
$\alpha_2'\alpha_2^{-1} \rtimes \nu^{1/2}\alpha_2 1_{\GSp(2)}$&IIIb \\
\hline
$\alpha_1'\alpha_2^{-1}= \nu$, $\alpha_2' \alpha_2^{-1} = \nu$ &
$\nu \rtimes \nu^{1/2}\alpha_2 1_{\GSp(2)}$&IIIb \\
\hline
$\alpha_1'\alpha_2^{-1}= \nu$, $\alpha_2' \alpha_2^{-1} = \nu^{-1}$ &
$\nu \rtimes \nu^{-1/2}\alpha_1 1_{\GSp(2)}$&IIIb \\
\hline
$\alpha_1'\alpha_2^{-1}= \nu^{-1}$, $\alpha_2' \alpha_2^{-1} \neq \nu^{\pm 1}$ &
$\alpha_1'\alpha_1^{-1} \rtimes \nu^{1/2}\alpha_1 1_{\GSp(2)}$&IIIb \\
\hline
$\alpha_1'\alpha_2^{-1}= \nu^{-1}$, $\alpha_2' \alpha_2^{-1} = \nu$&
$\nu\rtimes \nu^{-1/2}\alpha_1 1_{\GSp(2)}$& IIIb \\
\hline
$\alpha_1'\alpha_2^{-1}= \nu^{-1}$,  $\alpha_2' \alpha_2^{-1} = \nu^{-1}$&
$\nu\rtimes \nu^{1/2}\alpha_1 1_{\GSp(2)}$& IIIb \\
\hline
\multicolumn{3}{l}{\textbullet\qquad\qquad{$\pi_1 = \alpha_1 \times \alpha_2$\qquad $\pi_2 = \alpha 1_{\GL(2)}$ }}\\
\hline
$\alpha \alpha_2^{-1}  \neq \nu^{\pm 3/2}$ & $\alpha \alpha_2^{-1} 1_{\GL(2)} \rtimes \alpha_2$ & IIb\\
\hline
$\alpha \alpha_2^{-1} = \nu^{\pm 3/2}$ & $\alpha 1_{\GSp(4)}$ & IVd \\ 
\hline
\multicolumn{3}{l}{\textbullet\qquad\qquad{$\pi_1 = \alpha_1 \times \alpha_2$\qquad $\pi_2 = \alpha \St_{\GL(2)}$} }\\
\hline
$\alpha \alpha_2^{-1}  \neq \nu^{\pm 3/2}$ & $\alpha \alpha_2^{-1} \St_{\GL(2)} \rtimes \alpha_2$ & IIa\\
\hline
$\alpha \alpha_2^{-1} = \nu^{\pm 3/2}$ & $L(\nu^{3/2} \St_{\GL(2)}, \nu^{-3/2} \alpha)$ & IVc \\ 
\hline
\multicolumn{3}{l}{\textbullet\qquad\qquad{ $\pi_1 = \alpha_1 \times \alpha_2$\qquad $\pi_2$}  supercuspidal}\\
\hline
none & $\alpha_2^{-1} \pi_2\rtimes\alpha_2$ & X\\
\hline
\multicolumn{3}{l}{\textbullet\qquad\qquad{ $\pi_1 = \alpha 1_{\GL(2)}$ \qquad $\pi_2 =\beta 1_{\GL(2)}$}}\\
\hline
$\alpha \neq \beta$ & $L(\nu \alpha^{-1}\beta, \alpha^{-1} \beta \rtimes \nu^{-1/2} \alpha)$& Vd \\
\hline
$\alpha = \beta$ & $L(\nu, 1_{F^\times} \rtimes \nu^{-1/2} \alpha)$ & VId \\
\hline
\multicolumn{3}{l}{\textbullet\qquad\qquad{ $\pi_1 = \alpha 1_{\GL(2)}$ \qquad $\pi_2 =\beta \St_{\GL(2)}$}}\\
\hline
$\alpha \neq \beta$ & $L(\nu^{1/2} \alpha^{-1} \beta \St_{\GL(2)} , \nu^{-1/2} \alpha)$ & Vb \\
\hline
$\alpha =\beta$ & $L(\nu^{1/2} \St_{GL(2)}, \nu^{-1/2} \alpha)$& VIc \\
\hline
\multicolumn{3}{l}{\textbullet\qquad\qquad{$\pi_1 = \alpha 1_{\GL(2)}$\qquad $\pi_2$}  supercuspidal}\\
\hline
none & $L(\nu^{1/2}\alpha^{-1} \pi_2, \nu^{-1/2} \alpha)$ & XIb \\
\hline
\multicolumn{3}{l}{\textbullet\qquad\qquad{$\pi_1 = \alpha \St_{\GL(2)}$\qquad $\pi_2 = \beta \St_{\GL(2)}$}}\\
\hline
$\alpha \neq \beta$ & $\delta([\alpha^{-1} \beta, \nu \alpha^{-1} \beta ], \nu^{-1/2} \alpha)$, 
 supercuspidal& Va, Vb${}^*$\\
\hline
$\alpha=\beta$ & $\tau (S, \nu^{-1/2} \alpha)$, $\tau(T, \nu^{-1/2} \alpha)$ & VIa, VIb \\
\hline
\multicolumn{3}{l}{\textbullet\qquad\qquad{$\pi_1 = \alpha \St_{\GL(2)}$\qquad $\pi_2$} supercuspidal}\\
\hline
none&$\delta(\nu^{1/2}\alpha^{-1} \pi_2, \nu^{-1/2}\alpha)$, supercuspidal& XIa, XIb${}^*$ \\
\hline
\multicolumn{3}{l}{\textbullet\qquad\qquad$\pi_1$ supercuspidal\qquad $\pi_2$ supercuspidal}\\
\hline
$\pi_1 \ncong \pi_2$ &
 two element supercuspidal $L$-packet & XIIa${}^*$, XIIb${}^*$ \\
\hline
$\pi_1 \cong \pi_2$ & $\tau(S,\pi_1)$, $\tau(T,\pi_1)$ & VIIIa, VIIIb \\
\hline
\multicolumn{3}{l}{\textbullet\qquad\qquad
{$\pi_0 = \alpha_1 \times \alpha_2$}}\quad If $\alpha_1^\sigma \neq \alpha_1$, then $\xi\neq1$ is a character with $\xi^2=1$ and $\xi\pi(\alpha_1)\simeq\pi(\alpha_1)$.\\
\hline
$\alpha_1^\sigma \neq \alpha_1$, $\omega_{E/F}\eta^{-1} \alpha_2|_{F^\times} \neq \begin{cases}\xi\nu^{\pm 1}\\ 1 \end{cases}$& $ \omega_{E/F}\eta^{-1} \alpha_2|_{F^\times}  \rtimes \pi( \alpha_1)$ & VII\\
\hline
$\alpha_1^\sigma \neq \alpha_1$, $\omega_{E/F}\eta^{-1} \alpha_2|_{F^\times} =\xi\nu$&$L(\xi\nu,\pi(\alpha_1))$&IXb\\
\hline
$\alpha_1^\sigma \neq \alpha_1$, $\omega_{E/F}\eta^{-1} \alpha_2|_{F^\times} =\xi\nu^{-1}$ &$L(\xi\nu,\pi(\alpha_2))$&IXb\\
\hline
$\alpha_1^\sigma \neq \alpha_1$, $\omega_{E/F}\eta^{-1} \alpha_2|_{F^\times} =1$&$\tau(S,\pi(\alpha_2))$, $\tau(T,\pi(\alpha_2))$& VIIIa,VIIIb\\
\hline
$\alpha_1^\sigma=\alpha_1$, $\hat\alpha_1\circ\Norm^E_F=\alpha_1$&$\omega_{E/F}\eta^{-1}\hat\alpha_1^2\times\omega_{E/F}\rtimes\eta\hat\alpha_1^{-1}$&I\\
\hline
\multicolumn{3}{l}{\textbullet\qquad\qquad{$\pi_0=\alpha 1_{\GL(2)}$}}\\
\hline
$\alpha^\sigma\neq\alpha$& $L(\nu\omega_{E/F}\eta^{-1}\alpha|_{F^\times},\nu^{-1/2}\pi(\alpha))$&IXb\\
\hline
$\alpha^\sigma=\alpha$, $\hat\alpha\circ\Norm^E_F=\alpha$, $\hat\alpha^2=\eta$& $L(\nu\omega_{E/F},\omega_{E/F}\rtimes\nu^{-1/2}
\hat\alpha)$&Vd\\
\hline
$\alpha^\sigma=\alpha$, $\hat\alpha\circ\Norm^E_F=\alpha$, $\hat\alpha^2=\eta\omega_{E/F}$& $\omega_{E/F}\rtimes\hat\alpha 1_{\GSp(2)}$&IIIb\\
\hline
\multicolumn{3}{l}{\textbullet\qquad\qquad{$\pi_0=\alpha\St_{\GL(2)}$}}\\
\hline
$\alpha^\sigma\neq\alpha$&$\delta(\nu\omega_{E/F}\eta^{-1}\alpha|_{F^\times},\nu^{-1/2}\pi(\alpha))$&IXa\\
\hline
$\alpha^\sigma=\alpha$, $\hat\alpha\circ\Norm^E_F=\alpha$, $\hat\alpha^2=\eta$& $\begin{array}{c}\delta([\omega_{E/F},\nu\omega_{E/F}],\nu^{-1/2}\hat\alpha),\\
\text{\rm supercuspidal rep'n}\end{array}$&Va,Vb${}^*$\\
\hline
$\alpha^\sigma=\alpha$, $\hat\alpha\circ\Norm^E_F=\alpha$, $\hat\alpha^2=\eta\omega_{E/F}$& $\omega_{E/F}\rtimes\hat\alpha \St_{\GSp(2)}$&IIIa\\
\hline
\multicolumn{3}{l}{\textbullet\qquad\qquad$\pi_0$ supercuspidal}\\
\hline
$\pi_0^\sigma\cong\pi_0$, BC $(\hat \pi_0)=\pi_0$, $\omega_{\hat \pi_0} = \eta\omega_{E/F}$ & $\omega_{E/F} \rtimes \hat \pi_0$ & VII  \\
\hline
$\pi_0^\sigma\ncong\pi_0$&$\begin{array}{c} \mathrm{one}\ \mathrm{element}\ \mathrm{supercuspidal}\\ L\mathrm{-packet}\end{array}$ & XIII${}^*$\\
\hline
$\pi_0^\sigma\cong\pi_0$, BC$(\hat \pi_0)=\pi_0$, $\omega_{\hat \pi_0} = \eta$ & $\begin{array}{c} \mathrm{two}\ \mathrm{element}\ \mathrm{supercuspidal}\\ L\mathrm{-packet} \end{array}$& XIIa${}^*$, XIIb${}^*$ \\
\end{longtable}
\end{center}

{\bf Table 1 notation.}
If $\alpha$ is a non-Galois invariant character of $E^\times$, then $\pi(\alpha)$ is the associated supercuspidal representation of $\GL(2,F)$. 
If $\alpha$ is a Galois-invariant character of $E^\times$, then $\hat \alpha$ denotes a character of $F^\times$ such that $\hat \alpha \circ \Norm_F^E = \alpha$;  if $\pi_0$ is Galois-invariant, then   $\hat{\pi_0}$ denotes an irreducible, admissible representation of $\GL(2,F)$ such that the base change $\mathrm{BC}(\hat{\pi_0})$ of $\hat \pi_0$ is $\pi_0$.  The central character of a representation $\pi$ is denoted by  $\omega_\pi$.

\subsection*{Degree four invariants} 
We now specialize to the case in which $\pi_1$, $\pi_2$, and $\pi_0$ are tempered representations with trivial central character.  In this case, the $L$-packet $\Pi(\varphi)$ associated to $\varphi$, as described in the previous section, contains a unique generic representation $\pi$, which is also tempered and has trivial central character.  The following proposition tabulates the level $N_\pi$, the Atkin-Lehner eigenvalue $\varepsilon_\pi$, and the Hecke eigenvalues $\mu_\pi$ and $\lambda_\pi$ of the paramodular newform in $\pi$.  We also fix two additive characters, $\psi$ of $F$ with conductor $\OF$ and $\psi_E$ of $E$ with conductor $\OF_E$. For the purposes of stating the next proposition we define certain Euler-type factors. Let $N$ be a non-negative integer, let $\varepsilon = \pm 1$, and  $\lambda$ and $\mu$ be complex numbers. We define the factor $L(s,N,\varepsilon,\lambda,\mu)$ as follows. If $N=0$, then we define:
$$
L(s,0,\varepsilon,\lambda,\mu)^{-1}=
1-q^{-3/2}\lambda q^{-s} +(q^{-2}\mu+1+q^{-2})q^{-2s}-q^{-3/2} \lambda q^{-3s}+q^{-4s}
$$
If $N=1$, then we define
$$
L(s,1,\varepsilon,\lambda,\mu)^{-1}=
1-q^{-3/2}(\lambda+\varepsilon)q^{-s}+(q^{-2}\mu+1)q^{-2s}+\varepsilon q^{-1/2} q^{-3s}.
$$
If $N \geq 2$, then we define
$$
L(s,N,\varepsilon,\lambda,\mu)^{-1}=
1-q^{-3/2}\lambda q^{-s}+(q^{-2} \mu+1)q^{-2s}.
$$
\begin{proposition} 
\label{Hecketable}
Let the notation be as in Proposition \ref{Lpackettable}. Assume additionally that $\pi_1$, $\pi_2$ and $\pi_0$ are tempered, and that $\eta=1$. 
Then $\Pi(\varphi)$ contains a unique generic element $\pi$, which is tempered and  has trivial central character. 
Let $N_\pi$ be the paramodular level of $\pi$,  let $\varepsilon_\pi$ be the Atkin-Lehner eigenvalue of the newform in $\pi$, and let $\lambda_\pi$ and $\mu_\pi$ be the Hecke eigenvalues of the newform in $\pi$ for the Hecke operators $T_{0,1}$ and $T_{1,0}$ from 6.1 of \cite{RS1}. 
We have
\begin{equation}
\label{Npieq}
N_\pi = a ( \varphi)= 
\begin{cases}
a(\pi_1) + a(\pi_2) &\text{in the split case,}\\
2d(E/F)+f(E/F)a(\pi_0) &\text{in the non-split case},
\end{cases} 
\end{equation}
\begin{equation}
\label{epsiloneq}
\varepsilon_\pi = \varepsilon(1/2,\varphi, \psi,dx_{\psi})=
\begin{cases}
 \varepsilon(1/2,\pi_{1},\psi,dx_{\psi}) \varepsilon(1/2,\pi_{2},\psi,dx_{\psi})& \mbox{in the split case,}\\
 \varepsilon(1/2,\pi_{0},\psi_E,dx_{\psi_E})\omega_{E/F}(-1)&\mbox{in the non-split case,}
\end{cases}
\end{equation}
and
\begin{equation}
\label{Lequalitieseq}
 L(s,N_\pi, \varepsilon_\pi, \lambda_\pi,\mu_\pi) = L(s,\varphi) =
\begin{cases}
L(s,\pi_1)L(s,\pi_2)& \text{in the split case},\\
L(s,\pi_0) & \text{in the non-split case}.
\end{cases}
\end{equation}
Moreover, $N_\pi$, $\varepsilon_\pi$, $\lambda_\pi$ and $\mu_\pi$ are given by the following table.
Note that in the table, the subscript $\varpi$ indicates evaluation of the character at the uniformizer $\varpi$ of $F$.
\end{proposition}
\begin{center}
\renewcommand{\arraystretch}{1.5}
\setlength{\arraycolsep}{0.7mm}
\begin{longtable}{ccccc}
\caption{Hecke Eigenvalues}\\
\toprule
condition & $N_\pi$ & $\varepsilon_\pi$&$\lambda_\pi$&$\mu_\pi$ \\
\hline
\endfirsthead

\multicolumn{5}{c}{ \tablename\ \thetable{} -- continued from previous page}\\
\hline
condition & $N_\pi$ & $\varepsilon_\pi$&$\lambda_\pi$&$\mu_\pi$ \\
\endhead

\multicolumn{5}{c}{ Continued on next page}
\endfoot

\bottomrule
\endlastfoot
\multicolumn{5}{c}{\textbullet\qquad\qquad$\pi_1 = \alpha \times \alpha^{-1}$\quad  $\pi_2= \alpha' \times \alpha'^{-1}$\quad $\pi=\alpha \alpha' \times\alpha\alpha'^{-1} \rtimes \alpha^{-1}\qquad$\hfill I}\\
\hline
$\alpha, \alpha'$  unr &0&1 & 
$\begin{array}{l}
q^{3/2}(\alpha_\varpi+\alpha_\varpi^{-1}\\
\quad +\alpha'_\varpi+\alpha'_\varpi{}^{-1}) 
\end{array}$&
$\begin{array}{l}
q^2(1-q^{-2}+\\
(\alpha_\varpi+\alpha_\varpi^{-1})(\alpha'_\varpi+\alpha'_\varpi{}^{-1}))
\end{array}$\\
\hline
$\begin{array}{l}\alpha\ \mathrm{unr}\\\alpha'\ \mathrm{ram}\end{array}$& $2a(\alpha')$&$\alpha'(-1)$&$q^{3/2}(\alpha_\varpi+\alpha_\varpi^{-1})$ &0\\
\hline
$\begin{array}{l}\alpha\ \mathrm{ram}\\ \alpha'\ \mathrm{unr}\end{array}$&$2a(\alpha)$&$\alpha(-1)$&$q^{3/2}(\alpha'_\varpi+\alpha_\varpi'{}^{-1})$&0\\
\hline
$\alpha,\alpha'$ ram &$2(a(\alpha)+a(\alpha'))$&$\alpha(-1)\alpha'(-1)$&0&$-q^2$\\
\hline
\multicolumn{5}{c}{\textbullet\qquad\qquad$\pi_1=\alpha\times\alpha^{-1}$\quad $\pi_2=\alpha'\St_{\GL(2)}$\quad $\pi=\alpha\alpha'\St_{\GL(2)}\rtimes\alpha^{-1}$\hfill IIa}\\
\hline
$\alpha, \alpha'$  unr &1&$-\alpha'_\varpi$ & 
$\begin{array}{l}
q^{3/2}(\alpha_\varpi+\alpha_\varpi^{-1})\\
\quad +(q+1)\alpha'_\varpi 
\end{array}$&
$q^{3/2}(\alpha_\varpi\alpha'_\varpi+\alpha^{-1}_\varpi\alpha'_\varpi{}^{-1})$\\
\hline
$\begin{array}{l}\alpha\ \mathrm{unr}\\ \alpha'\ \mathrm{ram}\end{array}$& $2a(\alpha')$&$\alpha'(-1)$&$q^{3/2}(\alpha_\varpi+\alpha_\varpi^{-1})$ &0\\
\hline
$\begin{array}{l}\alpha\ \mathrm{ram}\\ \alpha'\ \mathrm{unr}\end{array}$   &$2a(\alpha)+1$&$-\alpha^{-1}(-1)\alpha'_\varpi$&$q\alpha'_\varpi$&$-q^2$\\
\hline
$\alpha,\alpha'$ ram&$2(a(\alpha)+a(\alpha'))$&$\alpha(-1)\alpha'(-1)$&0&$-q^2$\\
\hline
\multicolumn{5}{c}{\textbullet\qquad\qquad$\pi_1=\alpha \times \alpha^{-1}$\quad  $\pi_2$ supercuspidal\quad $\pi= \alpha \pi_2 \rtimes \alpha^{-1}$\hfill X}\\
\hline
$\alpha$ unr&$a(\pi_2)$ &$\varepsilon(\frac12,\pi_2)$&$q^{3/2}(\alpha_\varpi+\alpha_\varpi^{-1})$&0\\
\hline
$\alpha$ ram &$a(\pi_2)+2a(\alpha)$ &$\alpha(-1) \varepsilon(\frac12,\pi_2)$&0&$ -q^2$\\
\hline
\multicolumn{5}{c}{\textbullet\qquad\qquad$\pi_1 = \alpha \St_{\GL(2)}$\quad
$\pi_2 = \beta \St_{\GL(2)}$\quad $\alpha \neq \beta$\quad $\pi=
\delta([\alpha^{-1} \beta, \nu \alpha^{-1} \beta], \nu^{-1/2} \alpha)$\hfill Va}\\
\hline
$\alpha, \beta$ unr&2&-1&0&$-q^2-q$\\
\hline
$\begin{array}{c} \alpha\ \mathrm{unr} \\ \beta\ \mathrm{ram} \end{array}$ &$2a(\beta)+1$&$-\alpha_\varpi \beta(-1)$ &$\alpha_\varpi q$&$-q^2$\\
\hline
$\begin{array}{c} \alpha\ \mathrm{ram} \\ \beta\ \mathrm{unr} \end{array}$ &$2a(\alpha)+1$&$-\alpha(-1) \beta_\varpi$&$-\alpha_\varpi q$&$-q^2$\\
\hline
$\alpha, \beta$ ram &$2a(\alpha)+2a(\beta)$&$\alpha(-1)\beta(-1)$&0&$-q^2$\\
\hline
\multicolumn{5}{c}{\textbullet\qquad\qquad$\pi_1 = \alpha \St_{\GL(2)}$\quad
$\pi_2 = \beta \St_{\GL(2)}$\quad $\alpha = \beta$\quad $\pi=
\tau(S,\nu^{-1/2}\alpha)$\hfill VIa} \\
\hline
$\alpha$ unr&2&1&$2q\alpha_\varpi$&$-q(q-1)$\\
\hline
$\alpha$ ram&$4a(\alpha)$&1&0&$-q^2$\\
\hline
\multicolumn{5}{c}{\textbullet\qquad\qquad$\pi_1 = \alpha \St_{\GL(2)}$\quad
$\pi_2$ supercuspidal\quad $\pi=\delta( \nu^{1/2} \alpha^{-1} \pi_2, \nu^{-1/2} \alpha)$\hfill XIa} \\
\hline
$\alpha$ unr &$a(\pi_2)+1$ &$-\alpha_\varpi \varepsilon(\frac12,\pi_2)$&$q \alpha_\varpi$ &$-q^2$ \\
\hline
$\alpha$ ram & $a(\pi_2)+ 2 a(\alpha)$& $\alpha(-1) \varepsilon(\frac12,\pi_2)$ & 0 & $-q^2$ \\
\hline
\multicolumn{5}{c}{\textbullet\qquad\qquad $\pi_1$, $\pi_2$ supercuspidal\quad $\pi_1 \cong \pi_2$\quad $\pi = \tau(S,\pi_1)$\hfill VIIIa} \\
\hline
none&$2a(\pi_1)$ &1 & 0 & $-q^2$ \\
\hline
\multicolumn{5}{c}{\textbullet\qquad\qquad $\pi_1$, $\pi_2$ supercuspidal\quad $\pi_1 \ncong \pi_2$\quad $\pi$ supercuspidal\hfill XIIa${}^*$}\\
\hline
none& $a(\pi_1) + a(\pi_2)$ & $\varepsilon(\frac12, \pi_1) \varepsilon(\frac12, \pi_2)$ & 0 & $-q^2$ \\
\hline
\multicolumn{5}{c}{\textbullet\qquad\qquad$\pi_0=\alpha \times \alpha^{-1}$\quad
$\alpha^\sigma \neq \alpha$\quad $\alpha|_{F^\times} \neq \omega_{E/F}$\quad $\pi= \omega_{E/F}$ $\alpha|_{F^\times}^{-1} \rtimes \pi(\alpha^{-1})$\hfill VII} \\
\hline
none & $2d+2fa(\alpha)$ & $\omega_{E/F}(-1) \alpha(-1)$
& 0 & $-q^2$ \\
\hline
\multicolumn{5}{c}{\textbullet\qquad\qquad$\pi_0=\alpha \times \alpha^{-1}$\quad
$\alpha^\sigma \neq \alpha$\quad $\alpha|_{F^\times} = \omega_{E/F}$\quad $\pi= \tau(S,\pi(\alpha^{-1}))$\hfill VIIIa }\\
\hline
none & $2d+2fa(\alpha)$ & 1 & 0 & $-q^2$ \\
\hline
\multicolumn{5}{c}{\textbullet\qquad\qquad$\pi_0=\alpha \times \alpha^{-1}$\quad
$\alpha^\sigma = \alpha$\quad $\alpha = \hat \alpha \circ \Norm_F^E$\quad $\pi= \omega_{E/F} \hat \alpha^2 \times \omega_{E/F} \rtimes \hat \alpha^{-1}$\hfill I}\\
\hline$
\begin{array}{c}
 \alpha\ \mathrm{unr}\\
\omega_{E/F} \  \mathrm{unr}
\end{array}$
 &0&1& 0 &
$\begin{array}{c} -q^2  (\alpha_\varpi+ \alpha_\varpi^{-1})\\ -q^2-1 \end{array}$  \\
\hline
$\begin{array}{c}
\alpha\ \mathrm{unr}\\
\omega_{E/F}  \mathrm{ram}
\end{array}$
 & $2d$&$\omega_{E/F}(-1)$& $q^{3/2}( \hat{\alpha}_{\varpi} +  \hat{\alpha}_{\varpi}^{-1})$&0 \\
\hline
$ \alpha\ \mathrm{ram}$
 & $ 2d+2fa(\alpha)$ & $\omega_{E/F}(-1)$&0 & $-q^2$ \\
\hline
\multicolumn{5}{c}{\textbullet\qquad\qquad$\pi_0 = \alpha \St_{\GL(2)}$\quad $\alpha^\sigma \neq \alpha$\quad $\pi=\delta( \nu \omega_{E/F} \alpha|_{F^\times}$\quad $\nu^{-1/2} \pi(\alpha))$\hfill IXa} \\
\hline
none & $2d+2fa(\alpha)$ & $\omega_{E/F}(-1)\alpha (-1)$ & 0 & $-q^2$ \\
\hline
\multicolumn{5}{c}{\textbullet\qquad\qquad $\pi_0= \alpha \St_{\GL(2)}$\quad
$\alpha^\sigma = \alpha$\quad  $\alpha = \hat \alpha \circ \Norm_F^E$\quad $\hat \alpha^2 =1$\quad $\pi= \delta([\omega_{E/F},\nu \omega_{E/F}],\nu^{-1/2} \hat \alpha)$\hfill Va}\\
\hline
$\begin{array}{c} \alpha \ \mathrm{unr} \\ \omega_{E/F} \ \mathrm{unr} \end{array}$& 2 & -1 & 0 & $-q^2-q$  \\ 
\hline
$\begin{array}{c}  \alpha \ \mathrm{unr} \\ \omega_{E/F} \ \mathrm{ram} \end{array}$ & $2d+1 $& -  $\alpha_{\varpi_E} \omega_{E/F}(-1)$ &  $\alpha_{\varpi_E} q$& $-q^2$  \\ 
\hline
$ \alpha$ ram  & $2d+2fa(\alpha)$ &$\omega_{E/F}(-1)$& 0 & $-q^2$  \\ 
\hline
\multicolumn{5}{c}{\textbullet\qquad\qquad$\pi_0=\alpha\St_{\GL(2)}$\quad $\alpha^\sigma = \alpha$\quad $\alpha = \hat \alpha \circ \Norm_F^E$\quad $\hat \alpha^2 = \omega_{E/F}$\quad $\pi = \omega_{E/F} \rtimes \hat \alpha \St_{\GSp(2)}$\hfill IIIa}\\
\hline 
$\alpha$ unr&2 & 1 & $q(\hat{\alpha}_\varpi+\hat{\alpha}_\varpi^{-1})$&$-q^2+q$ \\
\hline
$ \alpha$ ram & $\begin{array}{c} 2d+2fa(\alpha)\\ =4a(\hat{\alpha})\end{array}$& 1 & 0 &$-q^2$ \\
\hline
\multicolumn{5}{c}{\textbullet\qquad\qquad$\pi_0\ \mathrm{supercuspidal}$\quad $\pi_0^\sigma \ncong \pi_0$\quad $\pi$ supercuspidal\hfill XIII${}^*$}\\
\hline
none& $2d +f a(\pi_0)$ & $\varepsilon(\frac12, \pi_0)\omega_{E/F}(-1)$ & 0 & $-q^2$ \\
\hline
\multicolumn{5}{c}{\textbullet\qquad\qquad $\pi_0$ supercuspidal\quad $\pi_0^\sigma \cong \pi_0$\quad  $\mathrm{BC}(\hat \pi_0) = \pi_0\quad \omega_{\hat \pi_0} = 1$\quad $\pi$ supercuspidal\hfill XIIa${}^*$}\\
\hline
none& $2d+f a(\pi_0)$ & $\varepsilon(\frac12, \pi_0)\omega_{E/F}(-1)$ & 0 & $-q^2$ \\
\hline
\multicolumn{5}{c}{\textbullet\qquad\qquad $\pi_0$ supercuspidal\quad $\pi_0^\sigma \cong \pi_0$\quad  $\mathrm{BC}(\hat \pi_0) = \pi_0$\quad $\omega_{\hat \pi_0} = \omega_{E/F}$\quad $\pi = \omega_{E/F} \rtimes \hat \pi_0$\hfill VII}\\
\hline
none& $\begin{array}{c}2d+fa(\pi_0)\\=2a(\hat \pi_0) \end{array}$& $\omega_{E/F}(-1)$ &0&$-q^2$ \\
\end{longtable}
\end{center}

\nll
{\bf Proof:} Let $\varphi=\varphi(\pi_1,\pi_2)$ or $\varphi(E,\pi_0,\eta)$ with $\pi_1$, $\pi_2$ and $\pi_0$ tempered,  $\omega_{\pi_1}=\omega_{ \pi_2}=1$, and $\eta=1$. Then using Proposition \ref{Lpackettable}, \cite{RS1} and the discussion preceding Proposition \ref{Lpackettable}, one can verify that all the elements of $\Pi(\varphi)$ are tempered with trivial central character, and  that exactly one element $\pi$ is generic.  Next,
the second equalities in (\ref{Npieq}) and (\ref{Lequalitieseq}) follow from $(a'1)$, $(a'2)$, $(L1)$ and $(L2)$ in \cite{Rohr}, along with the local Langlands correspondence for $\GL(2)$ (see, for example, \cite{Ku}). Similarly, the second equality in 
(\ref{epsiloneq}) follows from $(\epsilon'1)$ of \cite{Rohr} in the split case. Finally, assume that we are in the non-split case and let $\tilde{\psi}=\psi\circ\Trace^E_F$ . Then  we have
\begin{align*}
\varepsilon(\varphi,\psi,dx_{\psi}) 
&=
\varepsilon(\varphi_0,\tilde{\psi}, dx_{\tilde{\psi}}) 
\frac{\varepsilon(\Ind_{W_E'}^{W_F'} 1_E, \psi, dx_{\psi})^2}{\varepsilon(1_E, \tilde{\psi}, dx_{\tilde{\psi}})^2}\mbox{ (\cite{Rohr} } (\epsilon'2))\\
&=
\varepsilon(\varphi_0,\tilde{\psi}, dx_{\tilde{\psi}}) 
\frac{\varepsilon(1_F, \psi, dx_{\psi})^2\varepsilon(\omega_{E/F}, \psi, dx_{\psi})^2}{\varepsilon(1_E, \tilde{\psi}, dx_{\tilde{\psi}})^2}\\
&=
\varepsilon(\varphi_0,\psi_E, dx_{\psi_E}) 
\frac{\varepsilon(\omega_{E/F}, \psi, dx_{\psi})^2}{\varepsilon(1_E, \psi_E, dx_E)^2}\mbox{  (\cite{Rohr} 11 Prop.)}\\
&=\varepsilon(\varphi_0,\psi_E, dx_{\psi_E})\varepsilon(\omega_{E/F}, \psi, dx_{\psi})^2\\
&=\varepsilon(\varphi_0,\psi_E, dx_{\psi_E})\omega_{E/F}(-1)q^d \mbox{  (\cite{Rohr} 12 Lemma)}.
\end{align*}
Thus, at the center of the critical strip, we have again by \cite{Rohr} 11 Prop. (iii)
\begin{align*}
\varepsilon(1/2,\varphi,\psi,dx_{\psi}) &=\varepsilon(\varphi,\psi,dx_{\psi}) q^{-(2d+f a(\varphi_0))/2}\\
&=\varepsilon(\varphi_0,\psi_E, dx_{\psi_E})q_E^{-a(\varphi_0)/2}\omega_{E/F}(-1)\\
&=\varepsilon(1/2, \varphi_0,\psi_E, dx_{\psi_E})\omega_{E/F}(-1),
\end{align*}
as desired.
 \nl In the case that $\pi$ is non-supercuspidal the first equality in (\ref{Npieq}), (\ref{epsiloneq}) and (\ref{Lequalitieseq}) is known by \cite {RS1} Theorems 7.5.3,  7.5.9 and Corollary 7.5.5.  
\nl
Now consider the case that $\pi$ is a supercuspidal representation.  In this case all of the factors in (\ref{Lequalitieseq}) are 1.  Let $\epsilon(s,\pi)$ be the epsilon factor defined by the Novodvorsky zeta integrals of $\pi$ as discussed in Section \ref{supercuspidalepsilon}.  By Corollary 7.5.5 of \cite{RS1} we have that $\epsilon(s,\pi)=\epsilon_\pi q^{N_\pi(s-1/2)}$.  Further, as $\pi$ is supercuspidal, the epsilon factor and gamma factor defined by the Novodvorsky zeta integrals coincide.  By Proposition \ref{equalgammasplit} and Proposition \ref{equalgammainduced} we then have that
$$
\epsilon_\pi q^{N_\pi(s-1/2)}=\begin{cases}\gamma(s,\pi_1,\psi)\gamma(s,\pi_2,\psi)&\text{split}\\
                               \omega_{E/F}(-1)\gamma(s,\pi_0,\tilde{\psi})&\text{non-split}.
                              \end{cases}
$$
Finally, the proof is completed by noting that we have
\begin{align*}
\gamma(s,\pi_i,\psi)&=\epsilon(1/2,\pi_i,\psi,dx_\psi)q^{-a(\pi_i)(s-1/2)}\qquad i=1,2\quad \text{split,}\\
\gamma(s,\pi_0,\tilde{\psi})&=\epsilon(1/2,\pi_0,\psi_E,d\psi_E)q^{-(2d(E/F)+fa(\pi_0))(s-1/2)}\qquad\text{non-split.}
\end{align*}\qed

\section{Equality of gamma factors for supercuspidal representations}\label{supercuspidalepsilon}

The main result of this section is the calculation of the Novodvorsky gamma factors of the supercuspidal representations of type XIIa${}^*$ and XIII${}^*$. 
\nl
Let $\pi$ be a supercuspidal generic irreducible admissible representation of $\GSp(4,F)$ with trivial central character, and let $s\in\C$.  We say that $\pi$ admits an $s$-Bessel model if $\pi$ is isomorphic to a space of functions $B:\GSp(4,F)\rightarrow\C$ that satisfy
$$
B(\begin{bmatrix}1&&y_1&y_2\\&1&y_2&y_3\\&&1&\\&&&1\end{bmatrix}g)=\psi(y_2)B(g)
$$
and
$$
B(\begin{bmatrix}t_1&&&\\&t_2&&\\&&t_2&\\&&&t_1\end{bmatrix}g)=|t_2/t_1|^{s-1/2}B(g)
$$
for $y_1,y_2,y_3 \in F$, $t_1,t_2 \in F^\times$, and $g \in \GSp(4,F)$. 
If $\pi$ admits an $s$-Bessel model, then this model is unique (\cite{RS1} Prop. 2.5.7), and we denote it by $B_s(\pi)$.  
\nl 
Now, let $\pi$ be a representation of type XIIa${}^*$ so that there exist two nonisomorphic, supercuspidal, irreducible, admissible representations $\pi_1$ and $\pi_2$ of $\GL(2,F)$ with trivial central characters such that $\pi=\theta_X(\sigma)$ where $X$ is the four-dimensional hyperbolic quadratic space over $F$, and $\sigma$ is a representation of $\GO(X)$ constructed from $\pi_1$ and $\pi_2$.
Concretely, let $X=M_2(F)$, and  equip $X$ with the symmetric bilinear form defined by $(x,y)=\Trace(xy^*)/2$ where $\mat{a}{b}{c}{d}^*=\mat{d}{-b}{-c}{a}$.  Let $R:=\{(g,h)\in \GSp(4,F)\times\GO(X): \lambda(g)=\lambda(h)\}$, and let $\omega=\omega_\psi$ be the Weil representation of $R$ on the Schwartz space $\mathcal{S}(X^2)$ with respect to $\psi$ \cite{R}.  We have an exact sequence
$$
1\rightarrow F^\times\rightarrow \GL(2,F)\times\GL(2,F)\stackrel{\rho}{\rightarrow}\GSO(X)\rightarrow 1
$$
where $\rho(g_1,g_2)x=g_1xg_2^*$.  
Let $W(\pi_i,\psi)$ be the $\psi$-Whittaker model for $\pi_i$, so that $W_i \in W(\pi_i,\psi)$ transforms according to the formula
$$
W_i(\mat{1}{x}{}{1}g)=\psi(x)W_i(g)
$$
for $g\in\GL(2,F)$ and $x\in F$.
Let $x_1=\mat{0}{1}{0}{0}$ and $x_2=\mat{0}{0}{-2}{0}$ and let $H$ be the stabilizer of $x_1$ and $x_2$ in $\SO(X)$.
For $W_i\in W(\pi_i,\psi)$ and $\varphi\in S(X^2)$ we define
$$ 
B(g, \varphi, W_1, W_2, s):=\int\limits_{H\setminus\SO(X)}(\omega(g,hh^\prime)\varphi)(x_1,x_2)Z(s, \pi_1(h_1h_1^\prime)W_1)Z(s,\pi_2(h_2h_2^\prime)W_2)dh
$$
where $h=\rho(h_1,h_2)$, $h^\prime$ is any element of $\GSO(X)$ such that $\lambda(h^\prime)=\lambda(g)$ and
$$
Z(s,W_i)=\int_{F^\times}W_i(\mat{a}{}{}{1})|a|^{s-1/2}d^\times a,\quad W_i\in W(\pi_i,\psi).
$$
Note that as $\pi_i$ is supercuspidal, $Z(s,W_i)$ converges for all $s\in \C$ to a polynomial in $\C[q^{-s},q^{s}]$.  A similar statement is true for $B(g,\varphi, W_1, W_2, s)$.  One can prove that each of the functions $B(\cdot,\varphi,W_1, W_2,s)$ is contained in the $s$-Bessel model for $\pi$ and by extending linearly, we obtain a surjective map
\begin{equation*}
\beta_s:S(X^2)\otimes W(\pi_1,\psi)\otimes W(\pi_2,\psi)\rightarrow B_s(\pi)
\end{equation*}
with the property that for $g\in\GSp(4,F)$ and $h=\rho(h_1,h_2)\in\GSO(X)$ with $\lambda(h)=\lambda(g)$ we have
\begin{equation}\label{betatransform}
\beta_s(\omega(g,h)\varphi\otimes \pi_1(h_1)W_1\otimes \pi_2(h_2)W_2)=g\cdot\beta_s(\varphi\otimes W_1\otimes W_2).
\end{equation}
\nl 
On the other hand, for each $s\in\C$, there is another surjective  map 
\begin{equation*}
\beta_s^\prime:S(X^2)\otimes W(\pi_1,\psi)\otimes W(\pi_2,\psi)\rightarrow B_s(\pi)
\end{equation*}
with the analogous transformation property.  This map is constructed using the Weil representation and zeta integrals.  Let $c_1, c_2\in F^\times$ and let $W(\pi,\psi_{c_1,c_2})$ be the $\psi_{c_1,c_2}$-Whittaker model for $\pi$.  If $W\in W(\pi,\psi_{c_1,c_2})$, then 
\begin{equation}
\label{whiteq}
W(\begin{bmatrix}
1&x&*&*\\
&1&*&y\\
&&1&\\
&&-x&1  
  \end{bmatrix}g)=\psi(c_1x+c_2y)W(g)
\end{equation}
for $x,y\in F$ and $g\in \GSp(4,F)$.  
The map $\beta_s^\prime$ is defined to be the composition of $\GSp(4,F)$ maps
\begin{multline}\label{betacomp}
S(X^2)\otimes W(\pi_1,\psi)\otimes W(\pi_2,\psi)\stackrel{\mathrm{id}\otimes S_{1/2}}{\longrightarrow} S(X^2)\otimes W(\pi_1,\psi^{1/2})\otimes W(\pi_2,\psi^{-1/2})\\\stackrel{C}{\longrightarrow} W(\pi,\psi_{1/2,1/2})\stackrel{S_2}{\longrightarrow} W(\pi,\psi_{1,1})\stackrel{B_Z}{\longrightarrow} B_s(\pi).
\end{multline}
For the first map, we have
$$
S_{1/2}(W_1\otimes W_2)(g_1,g_2)=W_1(\mat{1/2}{}{}{1}g_1)\otimes W_2(\mat{-1/2}{}{}{1}g_2).
$$
For the second map, let $y_1=\mat{\ }{1}{}{}$ and $y_2=\mat{1}{}{}{1}$ and let $H^\prime$ be the stabilizer of $y_1$ and $y_2$ in $\SO(X)$.  Then the map $C$ is given by
$$
C(\varphi\otimes W_1\otimes W_2)(g)=\int\limits_{H^\prime\setminus\SO(X)}(\omega_\psi(g,hh^\prime)\varphi)(y_1,y_2)W_1(h_1h_1^\prime)W_2(h_2h_2^\prime)dh,
$$
where $h=\rho(h_1,h_2)$ and $h^\prime$ is any element of $\GSO(X)$ such that $\lambda(h^\prime)=\lambda(g)$.
The map $S_2$ is defined by the formula
$$
S_2(W)(g)=W(\begin{bmatrix}2&&&\\&1&&\\&&1/4&\\&&&1/2\end{bmatrix}g).
$$
To construct the final map, we recall that for $W\in W(\pi,\psi_{1,1})$, the zeta integral of $W$ is given by
$$
Z(s, W)=\int_{F^\times}\int_FW(\begin{bmatrix}a&&&\\&a&&\\&&1&\\&x&&1\end{bmatrix})|a|^{s-3/2}dxd^{\times}a.
$$
Since $\pi$ is supercuspidal, $Z(s,W)$ converges for all $s\in \C$ to a polynomial in $\C[q^{-s},q^s]$. We define the map $B_Z:W(\pi,\psi_{1,1})\rightarrow B_s(\pi)$ by 
$$
B_Z(W)(g):=Z(s,\pi(\begin{bmatrix}
                   1&&&\\
		   &&&1\\
		   &&1&\\
		   &-1&&
                  \end{bmatrix}g)W).
$$

\begin{lemma}\label{betabetaprimelemma}
There exists a constant $c\in \C^\times$ such that
$$
\beta_s(x)=c|2|^{s}\beta'_s(x)
$$
for all $x\in S(X^2)\otimes W(\pi_1,\psi)\otimes W(\pi_2,\psi)$ and for all $s\in\C$.
\end{lemma}

{\bf Proof:}  It follows from Theorem 1.8 of \cite{R} that for every $s\in \C$ the space of maps
$$
S(X^2)\otimes W(\pi_1,\psi)\otimes W(\pi_2,\psi)\rightarrow B_s(\pi)
$$
satisfying the transformation property in (\ref{betatransform}) is one dimensional.  Therefore, for every $s\in\C$ there exists a constant $c(s)\in\C^\times$ such that $\beta_s=c(s)\beta'_s$.  To compute $c(s)$, let $N_1$ be a positive integer and for $i=1,2$ let $W_i$ in the Whittaker model $W(\pi_i,\psi)$ of $\pi_i$ correspond to the
characteristic function $\chi_{1+\p^{N_1}}$ in the Kirillov model of $\pi_i$ with respect to $\psi$, so
that
$$
W_i(\mat{x}{}{}{1})=\chi_{1+\p^{N_1}}(x), \quad x \in F^\times.
$$
Choose $N_2>N_1$ such that
$$
\pi_i(\Gamma_{N_2}) W_i = W_i\qquad\text{for }i=1,2,
$$
where
$$
\Gamma_{N_2} = \{ \mat{a}{b}{c}{d} \in \GL(2,\OF): a,d \equiv 1\ (\p^{N_2}),\ b,c \equiv 0 \ (\p^{N_2}). \}.
$$
There is a homeomorphism
$$
H\setminus\SO(X)\stackrel{\sim}{\rightarrow}\SO(X)\cdot(x_1,x_2).
$$ 
Let $p$ be the projection
$p: \SO(X) \to H\backslash \SO(X)$. The set $p( \rho(\Gamma_{N_2} \times \Gamma_{N_2}) \cap \SO(X))$ is
an open neighborhood of $H \cdot 1$.  By applying the above  homeomorphism to this set one obtains an open neighborhood of $(x_1,x_2)$. Choose $N_3 > N_2$ such that this open neighborhood of $(x_1,x_2)$ contains
$$
\SO(X) (x_1,x_2) \cap (x_1 + \varpi^{N_3} \Mat(2,\OF), x_2 + \varpi^{N_3} \Mat(2,\OF)). 
$$
Let $\varphi =\varphi_1 \otimes \varphi_2 \in \mathcal{S}(X^2)$ be the characteristic function of $(x_1 + \varpi^{N_3} \Mat(2,\OF) ) \times ( x_2 + \varpi^{N_3} \Mat(2,\OF))$, where $\varphi_i$ is the characteristic function of $x_i + \varpi^{N_3} \Mat(2,\OF)$.  It follows that $\beta_s(\varphi\otimes W_1\otimes W_2)(1)$ is a non-zero constant $C_1$ independent of $s$.  A lengthy computation shows that $\beta'_s(\varphi\otimes W_1\otimes W_2)(1)=C_2|2|^{-s}$, where $C_2$ is a constant independent of $s$.  Thus, $c(s)=C_1C_2^{-1}|2|^s$.\qed
\nll

\begin{lemma}\label{betafunctionaleq}
For any $x\in S(X^2)\otimes W(\pi_1,\psi)\otimes W(\pi_2,\psi)$ and for any $g\in \GSp(4,F)$, there is a functional equation 
$$
\beta_{1-s}(x)(\begin{bmatrix}&1&&\\1&&&\\&&&1\\&&1&\end{bmatrix}g)=|2|^{1-2s}\gamma(s,\pi_1)\gamma(s,\pi_2)\beta_s(x)(g),
$$
where $\gamma(s, \pi_i)$ is the gamma factor for the $\GL(2,F)$ representation $\pi_i$.
\end{lemma}
\nll
{\bf Proof:} By (\ref{betatransform}), we may assume that $g=1$ and $x=\varphi\otimes W_1\otimes W_2$.  Note that
\begin{equation}\label{switchx1x2}
\rho(\mat{}{1/2}{-2}{},\mat{}{1}{-1}{}) x_1 = x_2, \quad \rho(\mat{}{1/2}{-2}{},\mat{}{1}{-1}{}) x_2 = x_1. 
\end{equation}
Then, 
\begin{align*}
&\beta_{1-s}(x)(\begin{bmatrix}&1&&\\1&&&\\&&&1\\&&1&\end{bmatrix})\\
&=\int_{H\backslash \SO(X)} \omega(1 ,\rho(h_1,h_2))\varphi (x_2,x_1) Z(1-s, \pi_1(h_1) W_1) Z(1-s, \pi_2(h_2) W_2)\, dh \\
&=\int_{H\backslash \SO(X)} \omega(1 ,\rho(h_1,h_2))\varphi (x_1,x_2)\\ &\qquad\qquad\qquad Z(1-s, \pi_1(\mat{}{1/2}{-2}{}^{-1}) \pi_1(h_1) W_1) Z(1-s, \pi_2(\mat{}{1}{-1}{}^{-1}) \pi_2(h_2) W_2)\, dh,
\end{align*}
by applying the identity (\ref{switchx1x2}).  For $i=1,2$, the zeta integral of $\pi_i$ satisfies a functional equation 
$$
Z(1-s,\pi_i(\mat{}{1}{-1}{})W)=\gamma(s,\pi_i)Z(s,W).
$$
This functional equation, together with the fact that these representations have trivial central character yields
$$
Z(1-s, \pi_1(\mat{}{1/2}{-2}{}^{-1})\pi_1(h_1) W_1)=|4|^{1-s-1/2}\gamma(s, \pi_1) Z(s, \pi_1(h_1) W_1), 
$$
and
$$
Z(1-s, \pi_2(\mat{}{1}{-1}{}^{-1}) \pi_2(h_2)W_2)=\gamma(s,\pi_2) Z(s,\pi_2(h_2) W_2). 
$$
Substituting these identities we obtain the lemma.\qed
\nl
\noindent We recall that the Novodvorsky zeta integrals for $\pi$ satisfy a functional equation.  In general if $\pi'$ is a generic irreducible admissible representation of $\GSp(4,F)$ with trivial central character, then there exists an element $\gamma(s,\pi')\in \C(q^{-s})$ such that 
$$
Z(1-s,\pi'(\begin{bmatrix}&&&1\\&&-1&\\&-1&&\\1&&&\end{bmatrix})W)=\gamma(s,\pi')Z(s,W)
$$
for all $W\in W(\pi',\psi_{1,1})$. See (\cite{RS1} Proposition 2.6.5).  Moreover the zeta integrals also define a local $L$-factor, $L(s,\pi')$, and epsilon factor 
$$
\epsilon(s,\pi')=\gamma(s,\pi')\frac{L(s,\pi')}{L(1-s,\pi')}.
$$
If $\pi'$ is supercuspidal, then $L(s,\pi')=1$ so that the epsilon and gamma factors coincide.
\begin{proposition}\label{equalgammasplit}
We have
$$
\gamma (s, \pi) = \gamma(s, \pi_1) \gamma(s, \pi_2). 
$$
\end{proposition}
\nll
\noindent {\bf Proof:} The proposition is proved by applying the previous lemmas to the functional equation for Novodvorsky zeta integrals.  
Let $x\in S(X^2)\otimes W(\pi_1,\psi)\otimes W(\pi_2,\psi)$.  Let $W^\prime \in W(\pi,\psi_{1,1})$ be the image of $x$ under the composition of the first three maps in (\ref{betacomp}).  The functional equation implies that $g\in \GSp(4,F)$ and for all $s\in\C$
\begin{align*}
Z(1-s, \pi(\begin{bmatrix}&&&1\\&&-1&\\&-1&&\\1&&&\end{bmatrix})\pi(\begin{bmatrix}1&&&\\&&&1\\&&1&\\&-1&&\end{bmatrix}g)W') = \gamma(s,\pi) Z(s, \pi(\begin{bmatrix}1&&&\\&&&1\\&&1&\\&-1&&\end{bmatrix}g)W'). 
\end{align*}
 Moreover, Lemma \ref{betabetaprimelemma} asserts that for all $g\in \GSp(4,F)$ and for all $s\in\C$,
\begin{equation}\label{betatozetaeq}
Z(s, \pi(\begin{bmatrix}1&&&\\&&&1\\&&1&\\&-1&&\end{bmatrix}g) W') = c^{-1}|2|^{-s} \beta_s(x)(g).
\end{equation}
Now the left hand side of the functional equation can be rewritten as
\begin{align*}
&Z(1-s, \pi(\begin{bmatrix}&&&1\\&&-1&\\&-1&&\\1&&&\end{bmatrix})\pi(\begin{bmatrix}1&&&\\&&&1\\&&1&\\&-1&&\end{bmatrix}g)W')\\
&= Z(1-s, \pi(\begin{bmatrix}1&&&\\&&&1\\&&1&\\&-1&&\end{bmatrix})\pi(\begin{bmatrix}&1&&\\1&&&\\&&&1\\&&1&\end{bmatrix} g) W') \\
&= c^{-1}|2|^{s-1} \beta_{1-s}(x)(\begin{bmatrix}&1&&\\1&&&\\&&&1\\&&1&\end{bmatrix} g).
\end{align*}
Substituting (\ref{betatozetaeq}) for the right hand side of the functional equation we obtain
$$
c^{-1}|2|^{s-1} \beta_{1-s}(x)(\begin{bmatrix}&1&&\\1&&&\\&&&1\\&&1&\end{bmatrix} g) = \gamma(s,\pi)c^{-1}|2|^{-s} \beta_s(x)(g).
$$
Applying Lemma \ref{betafunctionaleq}, we get:
\begin{align*}
c^{-1} |2|^{s-1} |2|^{1-2s}\gamma(s,\pi_1)\gamma(s,\pi_2)\beta_s(x)(g)
&=\gamma(s,\pi)c^{-1}|2|^{-s} \beta_s(x)(g)\\
\gamma(s,\pi_1)\gamma(s,\pi_2) \beta_s(x)(g)
&= \gamma(s,\pi) \beta_s(x)(g),
\end{align*}
for all $g\in \GSp(4,F)$ and for all $s\in \C$.  As $x$ and $g$ run through all possible values, the holomorphic functions in $s$, $\beta_s(x)(g)$, run through all possible zeta integrals of $\pi$ by (\ref{betatozetaeq}) and hence all functions of the form $P(q^{-s},q^s)$ where $P\in\C[X,Y]$ (\cite{RS1} Prop. 2.6.4; recall that $L(s,\pi)=1$).  Thus, we obtain the statement of the proposition. \qed

\nll
Now, suppose that $\pi$ is a representation of type XIII${}^*$. The shape of the argument is analogous to the previous case, but there are important differences which we will highlight.  In this case, there exists a quadratic extension $E=F(\sqrt{d_0})$ such that $\pi=\theta_{X_0}(\sigma_0)$, where $X_0$ is the four dimensional quadratic space having anisotropic component $(E,\Norm_F^E)$ and $\sigma_0$ is the supercuspidal, irreducible, admissible representation of $\GO(X_0)$ with trivial central character associated to  a supercuspidal, irreducible, admissible representation $\pi_0$ of $\GL(2,E)$ with trivial central character which is not Galois invariant.
Concretely, we take $X_0$ to be the subspace of $M_2(E)$ such that for all $ x\in X_0$, $\sigma(x)=x^*$, where $\mat{a}{b}{c}{d}^*=\mat{d}{-b}{-c}{a}$ and equip $X_0$ with a symmetric bilinear form $<x,y>=\Trace(xy^*)/2$.  Let $R:=\{(g,h)\in \GSp(4,F)\times\GO(X_0): \lambda(g)=\lambda(h)\}$, and let $\omega=\omega_\psi$ be the Weil representation of $R$ on the Schwartz space $\mathcal{S}(X_0^2)$ with respect to $\psi$.  We have an exact sequence
$$
1\rightarrow E^\times\rightarrow F^\times\times\GL(2,E)\stackrel{\rho_0}{\rightarrow}\GSO(X_0)\rightarrow 1
$$
where $\rho_0(t,g)x=t^{-1}gx\sigma(g)^*$ and the inclusion of $E^\times$ sends $z$ to $(N^E_F(z),z)$.  
Set $\tilde \psi = \psi \circ \Trace_{E/F}$. 
Let $W(\pi_0,\tilde{\psi})$ be the $\tilde{\psi}$-Whittaker model for $\pi_0$ such that $W_0 \in W(\pi_0,\tilde{\psi})$ transforms according to the formula
$$
W_0(\mat{1}{x}{}{1}g)=\tilde{\psi(x)}W_0(g)
$$
for $g\in\GL(2,E)$ and $x\in E$.
Let $x_1=\mat{0}{\sqrt{d_0}}{0}{0}$ and $x_2=\mat{0}{0}{-2/\sqrt{d_0}}{0}$ and let $H$ be the stabilizer of $x_1$ and $x_2$ in $\SO(X_0)$. For $g\in\GSp(4,F)^+$, $W_0\in W(\pi_0,\tilde{\psi})$ and $\varphi\in S(X_0^2)$ we define
$$ 
B(g, \varphi, W_0, s):=\int\limits_{H\setminus\SO(X)}(\omega(g,hh^\prime)\varphi)(x_1,x_2)Z(s, \pi_0(h_0h_0^\prime)W_0)dh
$$
where $h=\rho_0(t,h_0)$, $h^\prime$ is any element of $\GSO(X)$ such that $\lambda(h^\prime)=\lambda(g)$ and
$$
Z(s,W_0)=\int_{E^\times}W_0(\mat{a}{}{}{1})|a|_E^{s-1/2}d^\times a,\quad W_0\in W(\pi_0,\tilde{\psi}).
$$
Note that as $\pi_0$ is supercuspidal, $Z(s,W_0)$ converges for all $s\in \C$ to a polynomial in $q_E^{-s}$ and $q_E^{s}$.  A similar statement is true for $B(g,\varphi, W_0, s)$.  We extend $B(\cdot,\varphi,W_0,s)$ to $\GSp(4,F)$ via the formula
$$
B(g)=|\lambda(g)|^{-(s-1/2)}B(g_0)
$$
for 
$$
g_0=\begin{bmatrix}\lambda(g)^{-1}&&&\\&1&&\\&&1&\\&&&\lambda(g)^{-1}\end{bmatrix}g
$$
for $g \in \GSp(4,F)$. 
One can then prove that each of the functions $B(\cdot,\varphi,W_0,s)$ is contained in the $s$-Bessel model for $\pi$ and by extending linearly, we obtain a surjective map
\begin{equation*}
\beta_s:S(X^2)\otimes W(\pi_0,\tilde{\psi})\rightarrow B_s(\pi)
\end{equation*}
with the property that for $g\in\GSp(4,F)$ and $h=\rho_0(t,h_0)\in\GSO(X_0)$ with $\lambda(h)=\lambda(g)$ we have
\begin{equation}\label{betatransformnonsplit}
\beta_s(\omega(g,h)\varphi\otimes \pi_0(h_0)W_0=g\cdot\beta_s(\varphi\otimes W_0).
\end{equation}
\nl For each $s\in\C$, there is another surjective  map 
\begin{equation*}
\beta_s^\prime:S(X^2)\otimes W(\pi_0,\tilde{\psi})\rightarrow B_s(\pi)
\end{equation*}
with the analogous transformation property.  This map is constructed using the Weil representation and zeta integrals.  Let $c_1, c_2\in F^\times$ and let $W(\pi,\psi_{c_1,c_2})$ be the $\psi_{c_1,c_2}$-Whittaker model for $\pi$ as in (\ref{whiteq}). 
The map $\beta_s^\prime$ is defined to be the composition of $\GSp(4,F)$ maps
\begin{multline}\label{betacompnonsplit}
S(X^2)\otimes W(\pi_0,\tilde{\psi})\stackrel{\mathrm{id}\otimes S_{1/2\sqrt{d_0}}}{\longrightarrow} S(X^2)\otimes W(\pi_0,\tilde \psi^{1/2\sqrt{d_0}})\\\stackrel{C}{\longrightarrow} W(\pi,\psi_{1/2,1/2})\stackrel{S_2}{\longrightarrow} W(\pi,\psi_{1,1})\stackrel{B_Z}{\longrightarrow} B_s(\pi).
\end{multline}
For the first map, we have
$$
S_{1/2\sqrt{d_0}}(W_0)(g_0)=W_0(\mat{1/2\sqrt{d_0}}{}{}{1}g_0).
$$
For the second map, let $y_1=\mat{\ }{\sqrt{d_0}}{}{}$ and $y_2=\mat{1}{}{}{1}$ and let $H^\prime$ be the stabilize of $y_1$ and $y_2$ in $\SO(X)$.  Then the map $C$ is given by
$$
C(\varphi\otimes W_0)(g)=\int\limits_{H^\prime\setminus\SO(X)}(\omega_\psi(g,hh^\prime)\varphi)(y_1,y_2)W_0(h_0h_0^\prime)dh,
$$
where $g\in\GSp(4,F)^+$, $h=\rho(t,h_0)$ and $h^\prime$ is any element of $\GSO(X)$ such that $\lambda(h^\prime)=\lambda(g)$.  We extend the function $C(\varphi\otimes W_0)$ to all of $\GSp(4,F)$ by zero.  The maps $S_2$ and $B_Z$ are as in the previous case.

\begin{lemma}\label{betabetaprimelemmanonsplit}
There exists a constant $c\in \C^\times$ such that
$$
\beta_s(x)=c|2/d_0|^{s}\beta'_s(x)
$$
for all $x\in S(X^2)\otimes W(\pi_0,\tilde{\psi})$ and for all $s\in\C$.
\end{lemma}

{\bf Proof:}  The proof is analogous to the proof of Lemma \ref{betabetaprimelemma}. \qed

\nll

\begin{lemma}\label{betafunctionaleqnonsplit}
For any $x\in S(X^2)\otimes W(\pi_0,\tilde{\psi})$ and for any $g\in \GSp(4,F)$, there is a functional equation 
$$
\beta_{1-s}(x)(\begin{bmatrix}&1&&\\1&&&\\&&&1\\&&1&\end{bmatrix}g)=\omega_{E/F}(-1)|2/d_0|_E^{1/2-s}\gamma(s,\pi_0,\tilde{\psi})\beta_s(x)(g),
$$
where $\gamma(s, \pi_0, \tilde \psi)$ is the gamma factor for  $\pi_0$ with respect to $\tilde \psi$. 
\end{lemma}
\nll
{\bf Proof:} The proof is analogous to the proof of Lemma \ref{betafunctionaleq}.\qed

\begin{proposition}\label{equalgammainduced}
We have
$$
\gamma (s, \pi) = \omega_{E/F}(-1)\gamma(s, \pi_0,\tilde{\psi}). 
$$
\end{proposition}
\nll
\noindent {\bf Proof:} 
The proof is analogous to the proof of Proposition \ref{equalgammasplit}.\qed

\end{document}